\documentclass[12pt]{article}
\usepackage[T1]{fontenc}
\usepackage{amsmath}
\usepackage{amssymb}
\usepackage{latexsym}
\usepackage{textcomp}
\usepackage{diagrams}
\diagramstyle{h=21pt,w=21pt,scriptlabels,midshaft,PostScript=dvips,nohug,shortfall=4pt}

\usepackage[noBBpl]{mathpazo}
\renewcommand{\texttt}[1]{{\fontfamily{pcr}\fontseries{m}\fontshape{n}%
\selectfont #1}}

\newcommand{\kat}[1]{\text{\textbf{\textsl{#1}}}}
\renewcommand{\ldots}{\relax\ifmmode\ldotp\ldotp\ldotp\else$\m@th\ldotp\ldotp\ldotp\ $\fi}
\providecommand{\qed}{\hspace*{\fill}\nolinebreak[1]\hspace*{\fill} $\Box$}
\newcommand{\tensor}	{\otimes}
\newcommand{\ground}{\Bbbk}
\newcommand{\isopil}{\stackrel{\raisebox{0.1ex}[0ex][0ex]{\(\sim\)}}%
  {\raisebox{-0.15ex}[0.28ex]{\(\rightarrow\)}}}
\newcommand{\df}{\: {\raisebox{0.255ex}{\normalfont\scriptsize :\!\!}}=}

\newcommand{\wtil}{\widetilde}

\newcommand{\hovedfont}{\normalfont\bfseries}
\usepackage{theorem}
\theorembodyfont{\normalfont\itshape}
\theoremheaderfont{\hovedfont}
	\theoremstyle{change}
\newtheorem{lemma}{Lemma.}[section]
\newtheorem{prop}[lemma]{Proposition.}

\newtheorem{cor}[lemma]{Corollary.}
	\theorembodyfont{\normalfont}

\newtheorem{BM}[lemma]{Remark.}

\newtheorem{taller}[lemma]{$\!\!$}
\newenvironment{blanko}[1]%
{\begin{taller}{\hovedfont #1}\normalfont}%
{\end{taller}}
\newenvironment{deff}%
{%
\begin{list}{\em Definition. }%
{\setlength{\labelsep}{0mm}\setlength{\leftmargin}{0mm}%
\setlength{\labelwidth}{0mm}\setlength{\listparindent}{\parindent}%
\setlength{\parsep}{\parskip}\setlength{\partopsep}{0mm}}%
\item%
}%
{%
\end{list}%
}
\newenvironment{dem}%
{%
\begin{list}{\em Proof. }%
{\setlength{\labelsep}{0mm}\setlength{\leftmargin}{0mm}%
\setlength{\labelwidth}{0mm}\setlength{\listparindent}{\parindent}%
\setlength{\parsep}{\parskip}\setlength{\partopsep}{0mm}}%
\item%
}%
{%
\qed\end{list}%
}
\newenvironment{dem*}[1]%
{%
\begin{list}{\em #1 }%
{\setlength{\labelsep}{0mm}\setlength{\leftmargin}{0mm}%
\setlength{\labelwidth}{0mm}\setlength{\listparindent}{\parindent}%
\setlength{\parsep}{\parskip}\setlength{\partopsep}{0mm}}%
\item%
}%
{%
\qed\end{list}%
}
{%
\begin{list}{\em Proof. }%
{\setlength{\labelsep}{0mm}\setlength{\leftmargin}{0mm}%
\setlength{\labelwidth}{0mm}\setlength{\listparindent}{\parindent}%
\setlength{\parsep}{\parskip}\setlength{\partopsep}{0mm}}%
\item%
}%
{%
\qed\end{list}%
}
\newenvironment{bevis*}[1]%
{%
\begin{list}{\em #1 }%
{\setlength{\labelsep}{0mm}\setlength{\leftmargin}{0mm}%
\setlength{\labelwidth}{0mm}\setlength{\listparindent}{\parindent}%
\setlength{\parsep}{\parskip}\setlength{\partopsep}{0mm}}%
\item%
}%
{%
\qed\end{list}%
}

\newenvironment{blanko*}[1]%
{%
\begin{list}{\bf {#1} }%
{\setlength{\labelsep}{0mm}\setlength{\leftmargin}{0mm}%
\setlength{\labelwidth}{0mm}\setlength{\listparindent}{\parindent}%
\setlength{\parsep}{\parskip}\setlength{\partopsep}{0mm}}%
\item%
}%
{%
\end{list}%
}

\newcounter{dummycounter}
\newenvironment{punkt-a}%
{%
	\begin{list}%
	{(\alph{dummycounter})\hfill}%
	{\usecounter{dummycounter}%
	\setlength{\itemsep}{0em}\setlength{\parsep}{0em}\setlength{\topsep}{0em}%
	\setlength{\itemindent}{0em}\setlength{\labelwidth}{1.5em}
	\setlength{\labelsep}{0.3em}\setlength{\leftmargin}{1.8em}}%
}%
{\end{list}}

\makeatletter
\renewcommand{\ps@headings}
	{\setlength{\headheight}{16pt}%
	 \setlength{\headsep}{12pt}%
	 \renewcommand{\@oddhead}{\parbox{\textwidth}{%
			\small
			\texttt{Joachim Kock: Units \ \ \ \lastUpdate{2006-06-08 16:46}
			\hfill [\thepage/\pageref{lastpage}]}
			\\ \rule[8pt]{\textwidth}{0.3pt}}%
	 }
	\renewcommand{\@oddfoot}{}
	\renewcommand{\@evenfoot}{}%
}
\makeatother

\usepackage{mathrsfs}
\providecommand{\lastUpdate}[1]{#1}
\providecommand{\kat}[1]{\text{\textbf{\textsl{#1}}}}

\newcommand{\Vect}{\kat{Vect}}
\newcommand{\Nat}{\operatorname{Nat}}
\newcommand{\id}{\operatorname{id}}
\newcommand{\CC}{\mathscr{C}}
\newcommand{\DD}{\mathscr{D}}
\newcommand{\UU}{\mathscr{U}}

\begin{document}

% \renewcommand{\baselinestretch}{1.45}
%     \normalsize

\pagestyle{headings}

\def\vspec#1{\special{ps:#1}}%  passes #1 verbatim to the output
\def\rotstart#1{\vspec{gsave currentpoint currentpoint translate
	#1 neg exch neg exch translate}}% #1 can be any origin-fixing transformation
\def\rotfinish{\vspec{currentpoint grestore moveto}}% gets back in synch
	
\def\psrotate#1#2{\rotstart{#1 rotate}#2\rotfinish}

\def\psr#1{\rotstart{36.87 rotate}\hbox to0pt {\vsize 0pt \hss\(\scriptstyle #1\)\hss}\rotfinish}

% Scaling factor for finetuning globs.  Should always be 1.
\def\scaleFactor{10} % MUST BE AN INTEGER

\newcommand{\dropglob}[1]{%
\setlength{\unitlength}{0.003\DiagramCellWidth}
\multiply \unitlength by \scaleFactor
\begin{picture}(0,0)(0,0)
\qbezier(-28,-4)(0,-18)(28,-4)
\put(0,-14){\makebox(0,0)[t]{$\scriptstyle {#1}$}}
\put(28.6,-3.7){\vector(2,1){0}}
\end{picture}
}

\newcommand{\topglob}[1]{%
\setlength{\unitlength}{0.003\DiagramCellWidth}
\multiply \unitlength by \scaleFactor
\begin{picture}(0,0)(0,0)
\qbezier(-28,11)(0,25)(28,11)
\put(0,21){\makebox(0,0)[b]{$\scriptstyle {#1}$}}
\put(28.6,10.7){\vector(2,-1){0}}
\end{picture}
}

% The first argument is the lift (measured in pts)
\newcommand{\lift}[2]{%
\setlength{\unitlength}{1pt}
\begin{picture}(0,0)(0,0)
\put(0,{#1}){\makebox(0,0)[b]{${#2}$}}
\end{picture}
}

%%% OVERLINES %%%
\newlength{\ollength}

\newcommand{\ov}{\overline}

\newcommand{\OV}[1]{%
  \settowidth{\ollength}{\ensuremath{#1}}
  \addtolength{\ollength}{-2.5pt}
  \overset{\rule[-0.6pt]{1.5pt}{0pt}\rule[-0.6pt]{\ollength}{0.5pt}}{#1}%
}
\newcommand{\OVscript}[1]{%
  \settowidth{\ollength}{\ensuremath{#1}}
  \addtolength{\ollength}{-2.0pt}
  \setlength{\ollength}{0.75\ollength}
  \overset{\rule[-0.5pt]{1pt}{0pt}\rule[-0.5pt]{\ollength}{0.5pt}}{#1}%
}
\newcommand{\OVscriptscript}[1]{%
  \settowidth{\ollength}{\ensuremath{#1}}
  \addtolength{\ollength}{-1.5pt}
  \setlength{\ollength}{0.5\ollength}
  \overset{\rule[-0.4pt]{0.5pt}{0pt}\rule[-0.4pt]{\ollength}{0.5pt}}{#1}%
}
% We use a slightly longer overline for I:
\newcommand{\stregI}{%
  \settowidth{\ollength}{\ensuremath{I}}
  \addtolength{\ollength}{-2.0pt}
  \overset{\rule[-0.6pt]{1.5pt}{0pt}\rule[-0.6pt]{\ollength}{0.5pt}}{I}%
}
\newcommand{\stregIscript}{%
  \settowidth{\ollength}{\ensuremath{I}}
  \addtolength{\ollength}{-1.5pt}
  \setlength{\ollength}{0.75\ollength}
  \overset{\rule[-0.5pt]{1pt}{0pt}\rule[-0.5pt]{\ollength}{0.5pt}}{I}%
}
\newcommand{\stregIscriptscript}{%
  \settowidth{\ollength}{\ensuremath{I}}
  \addtolength{\ollength}{-1.0pt}
  \setlength{\ollength}{0.5\ollength}
  \overset{\rule[-0.4pt]{0.5pt}{0pt}\rule[-0.4pt]{\ollength}{0.5pt}}{I}%
}

\newcommand{\cel}[1]{\ensuremath{\mathsf{#1}}}

\newcommand{\LR}{LR}

\vspace*{24pt}

\begin{center}

  {\Large Elementary remarks on units in monoidal categories}

  \bigskip

  {\sc Joachim Kock}

\end{center}

\begin{abstract}
  We explore an alternative definition of unit in a monoidal category
  originally due to Saavedra: a Saavedra unit is a cancellative
  idempotent (in a $1$-categorical sense).  This notion is more
  economical than the usual notion in terms of left-right constraints,
  and is motivated by higher category theory.  To start, we describe
  the semi-monoidal category of all possible unit structures on a
  given semi-monoidal category and observe that it is contractible (if
  nonempty).  Then we prove that the two notions of units are
  equivalent in a strong functorial sense.  Next, it is shown that the
  unit compatibility condition for a (strong) monoidal functor is
  precisely the condition for the functor to lift to the categories of
  units, and it is explained how the notion of Saavedra unit naturally
  leads to the equivalent non-algebraic notion of fair monoidal
  category, where the contractible multitude of units is considered as
  a whole instead of choosing one unit.  To finish, the lax version of
  the unit comparison is considered.  The paper is self-contained.
  All arguments are elementary, some of them of a certain beauty.
\end{abstract}


\noindent
{\em Subject classification (MSC2000)}: 18D10

\normalsize

\makeatletter
\renewcommand{\section}{\@startsection {section}{1}{\z@}%
{-3.5ex \@plus -1ex \@minus -.2ex}%
{2.3ex \@plus.2ex}%
{\normalfont\large\bfseries}}
\renewcommand{\subsection}{\@startsection {subsection}{1}{\z@}%
{-3.5ex \@plus -1ex \@minus -.2ex}%
{2.3ex \@plus.2ex}%
{\normalfont\normalsize\bfseries}}

\setcounter{secnumdepth}{1}

\makeatother

%%%%%%%%%%%%%%%%%%%%%%%%%%%%%%%%%%%%%%%%%%%%%%%%%%
\section*{Introduction}
%%%%%%%%%%%%%%%%%%%%%%%%%%%%%%%%%%%%%%%%%%%%%%%%%%

\begin{blanko*}{Monoidal categories.}
  Monoidal categories are everywhere in mathematics, and serve among 
  other things as carrier for virtually all algebraic structures.
%   
%   
%   : every time you
%   have a category with products or sums or tensor products or direct
%   sums or disjoint unions or inter\-sections, etc.~--- whenever these
%   composition laws admit a neutral object.  It is the $1$-categorical
%   version of the notion of monoid, perhaps the most fundamental
%   algebraic structure in mathematics.
% 
  Monoidal categories are also the simplest example of higher
  categories, being the same thing as bicategories with only one
  object, just as a monoid can be seen as a category with only one
  object.  
%   Monoidal categories are quite easy to understand, but
%   higher categories are in general very hard.
  With the rapidly
  growing importance of higher category theory, it is interesting to
  revisit even the most basic theory of monoidal categories, to test
  new viewpoints and experiment with new formulations (cf.~also
  Chapter 3 in Leinster's book~\cite{Leinster-0305049}).
\end{blanko*}

\begin{blanko*}{Units.}
  This note analyses the notion of unit in monoidal categories.  Units
  enjoy a mixed reputation: in some monoidal categories, the unit
  appears to be an insignificant part of the structure and is often
  swept under the carpet; in other cases (like in linear logic, in the
  theory of operads, or in higher categories), the proper treatment of
  units can be highly non-trivial, suggesting that we have not yet 
  fully understood the nature of units.
  
  While the axioms for the multiplication law are rather well
  understood, and fit into a big pattern continuing in higher
  dimension (the geometrical insight provided by the Stasheff
  associahedra \cite{Stasheff:1963}), the unit axioms are subtler,
  and so far there seems to be no geometrical `explanation' of them.
  This delicacy is perhaps also reflected historically.
  The first finite list of axioms for a monoidal category was given by
  Mac Lane~\cite{MacLane:naturalAssociativity} in 1963, including one
  axiom for associativity (the pentagon equation) and four axioms for
  the unit with its left and right constraints.
  Shortly after, it was shown by Kelly~\cite{Kelly:MacLanesCoherence}
  that one of these four axioms for units in fact implies the three
  others.  His proof constitutes nowadays the first three lemmas in
  many treatments of monoidal categories, while other sources continue
  to employ the redundant axiomatics.  It is less well-known that
  conversely the three other unit axioms imply Kelly's axiom; this was
  observed by Saavedra~\cite{Saavedra} in 1972.
  
  However, mere rearrangement of the axioms imposed on the unit
  structure is not the crux of the matter.  The structure itself must
  be analysed.  As it turns out, the classical notion of unit is
  overstructured.  Saavedra~\cite{Saavedra} seems to have been the
  first to notice this: he observed that it is possible to express the
  notion of unit in monoidal categories without even mentioning the
  left and right constraints: a {\em Saavedra unit} (cf.~\ref{unit1})
  is an object $I$ equipped with an isomorphism $\alpha:I\tensor I
  \isopil I$, and having the property that tensoring with $I$ from
  either side is an equivalence of categories (in short, we say $I$ is
  a cancellable idempotent).  Saavedra observed that this notion is
  equivalent to the classical notion (although his proof has a gap, as
  far as I can see, cf.~Remark~\ref{error} below), but he did not
  pursue the investigation further --- he did not even consider
  monoidal functors in this viewpoint.

  The present note exploits the notion of Saavedra unit systematically
  to throw light also on the classical notion of unit and exhibit its
  redundancy.  Many mathematicians experience this redundancy even at
  a na\"\i ve level, writing for example `the right constraint is
  treated similarly', without realising that these phenomena can be
  distilled into precise results.  At a deeper level, one may expect 
  that the more economical notion of unit can help aliviate the 
  nastiness of units felt at times, since in general it is cheaper to 
  check a property than to provide a structure.
\end{blanko*}

\begin{blanko*}{Terminology.}
  Saavedra~\cite{Saavedra} used the term {\em reduced unit} for these
  cancellable idempotents, since they are less structured than the
  units-with-left-and-right-constraints.  Another option would be {\em
  absolute unit}, referring to the fact that the notion makes sense
  prior to any associativity constraints, so for example it makes
  sense to fix the Saavedra unit and vary the associators.  In
  contrast, the units axioms given in terms of left and right
  constraints make sense only relative to a fixed associativity
  constraint.  But in any case it seems unfortunate to differentiate
  this alternative notion of unit from the usual notion by means of an
  extra adjective, since the former are not a special kind of the
  latter -- the two notions are equivalent.  From a purist's
  viewpoint, the cancellable-idempotent notion is what really should
  be called {\em unit}, the richer structure of left and right
  constraints being something derived, as explained in
  Section~\ref{Sec:Saavedra}.  In this note, a temporary terminology
  is adopted where both notions of unit carry an extra attribute: the
  notion of unit in terms of left-right constraints is called {\em
  \LR{} unit}, while the cancellable-idempotent notion is called {\em
  Saavedra unit}.
\end{blanko*}

\begin{blanko*}{Overview of results.}
  The material is organised as follows.  In
  Section~\ref{Sec:classical}, after quickly reproducing Kelly's unit
  argument, we describe the category $U(\CC)$ of all possible \LR{}
  unit structures on a semi-monoidal category $\CC$ and show that it
  is contractible (if non-empty) (\ref{clasUcontract}).  We observe
  that $U(\CC)$ has a canonical semi-monoidal structure
  (\ref{tensor-unit}).

  In Section~\ref{Sec:Saavedra} we introduce the notion of Saavedra
  unit, and show how Saavedra units are canonically \LR{} units and
  vice versa (\ref{bij}).  We also define morphisms of Saavedra units,
  and show that the category of \LR{} units is isomorphic to the
  category of Saavedra units (\ref{unit-morphisms}).

  In Section~\ref{Sec:functoriality} we study (strong) monoidal
  functors, and show that compatibility with \LR{} units implies
  compatibility with Saavedra units and vice-versa
  (\ref{compatibilities}), and more precisely (\ref{2-iso}): there is
  an isomorphism between the $2$-category of monoidal categories with
  \LR{} units (and strong monoidal functors and monoidal natural
  transformations) and the $2$-category of monoidal categories with
  Saavedra units (and strong monoidal functors and monoidal natural
  transformations).  Two corollaries are worth mentioning: first
  (\ref{functor-redundance}), a strong monoidal functor is compatible
  with left constraints if and only if it is compatible with right
  constraints, and in fact this compatibility can be measured on $I$
  alone!  Second, a multiplicative functor is monoidal if and only if
  the image of a unit is again a unit (\ref{image-of-unit}).  This
  statement does not even make sense for \LR{} units (since the image
  of an \LR{} unit does not have enough structure to make sense of the
  question whether it is a unit again).  Finally we prove a rather
  technical result (\ref{extend}), needed in Section~\ref{Sec:gentle}:
  a unit compatibility on a strong monoidal functor $\CC\to\DD$ is
  equivalent to a lift to the categories of units $U(\CC)\to U(\DD)$.

  Section~\ref{Sec:monoids} is a short interlude on monoids, placed
  here in order to motivate the next notion: in
  Section~\ref{Sec:gentle} we introduce a relativisation of the notion
  of cancellable object, here called {\em gentle functors}: these are
  functors $\UU\to (\CC,\tensor)$ such that the two composites
  $\UU\times\CC\to\CC\times\CC\stackrel{\tensor}{\to} \CC
  \stackrel{\tensor}{\leftarrow}\CC\times\CC\leftarrow\CC\times\UU$
  are fully faithful.  It is easy to see that a Saavedra unit in $\CC$ is the
  same thing as a strong multiplicative functor $* \to \CC$ which is
  furthermore gentle.  Following this idea, we come to the notion of
  fair monoidal category: it is a gentle (strict) multiplicative
  functor $\UU\to \CC$ with $\UU$ contractible.  Here $\UU$ is thought
  of as the category of all units in $\CC$.  It is shown (\ref{fair}) that this
  notion of monoidal category is equivalent to the classical notion
  (as claimed in \cite{Kock:fair}).

  Finally, in Section~\ref{Sec:lax}, we generalise some of the results
  about strong monoidal functors to lax monoidal functors.  The
  Saavedra-unit compatibility is a bit less elegant to express in this
  case, but again we get an isomorphism of $2$-categories
  (\ref{lax-prop}).  As a particular case we get an isomorphism of the
  categories of monoids in the \LR{} and Saavedra-unit viewpoints,
  proving an assertion left open in Section~\ref{Sec:monoids}.
\end{blanko*}
 
\begin{blanko*}{Generalisations and outlook.}
  Throughout we assume the tensor product to be strict.  This is just
  for convenience: every argument can be carried over to the
  non-strict case, simply by inserting associator isomorphisms as
  needed.  That complication would not seem to illuminate anything
  concerning units.
  
  The notion of Saavedra unit has an obvious many-object version
  yielding an alternative notion of identity arrow in a bicategory.
  The notion of fair monoidal category described in
  Section~\ref{Sec:gentle} is just the one-object version of the
  notion of fair $2$-category, which has generalisations to fair
  $n$-categories~\cite{Kock:fair}.

  The remarks compiled in this note are a by-product of a more general
  investigation of weak units and weak identity arrows in higher
  categories~\cite{Kock:fair}, \cite{Joyal-Kock:coherence},
  \cite{Joyal-Kock:traintracks}, but as it turns out, even the
  $1$-dimensional case contains some surprises, and I found it
  worthwhile to write it down separately and explicitly, since I think
  it deserves a broader audience.
% 
%   and serves as a home for many easy arguments that did
%   not fit into the more advanced papers \cite{Kock:fair},
%   \cite{Joyal-Kock:coherence}, \cite{Joyal-Kock:traintracks}, but
%   which nevertheless are too cute to keep secret. 
  The notion of Saavedra unit dropped out of the theory of fair
  categories~\cite{Kock:fair}.  I am thankful to Georges
  Maltsiniotis for pointing out that the viewpoint goes back to
  Saavedra~\cite{Saavedra}.  In this note, reversing my own path into
  the subject, it is shown how the basic notion of Saavedra unit leads
  to the notion of fair (monoidal) category.  The two papers
  \cite{Joyal-Kock:coherence} and \cite{Joyal-Kock:traintracks} joint
  with Andr\'e Joyal deal with units in monoidal $2$-categories.  In
  \cite{Joyal-Kock:traintracks} we use the $2$-dimensional notion of
  Saavedra unit to prove a version of Simpson's weak-unit
  conjecture~\cite{Simpson:9810} in dimension $3$: strict
  $3$-groupoids with weak units can model all $1$-connected homotopy
  $3$-types.  The relevance of the Saavedra-unit viewpoint in
  higher-dimensional category theory was first suggested by
  Simpson~\cite{Simpson:9810}.  The main advantage is its economical
  nature, and in particular it is important that the very notion of
  unit is expressed in terms of a multiplication map, already a
  central concept in the whole theory.
\end{blanko*}

%   As
% a consequence of this, the microcosm yoga of monoids in monoidal
% categories which in turn are themselves monoids in $\Cat$ becomes even
% more accentuated, since units can be seens as certain semi-monoids in
% the dimension below: quoting from \ref{fair-monoid} and referring to
% Sections~\ref{Sec:monoids}~and~\ref{Sec:gentle} for the notion of
% gentle: `A fair monoidal category is a gentle semi-monoidal functor
% $\eta:U \to C$ with $U$ contractible.  A monoid in here is a gentle
% semi-monoid homomorphism $\eta(I) \to M$, where $I$ is a semi-monoid
% in $U$ and $M$ is a semi-monoid in $C$.'

\begin{blanko*}{Acknowledgments.}
  Part of this work was carried out while I was a postdoc at the
  Universit\'e du Qu\'ebec \`a Montr\'eal, supported by a CIRGET
  grant.  I wish to thank everybody at the UQ\`AM, and Andr\'e Joyal
  in particular, for a wonderful year in Montr\'eal.  Presently I am
  supported by a Ram\'on y Cajal fellowship from the Spanish Ministry
  of Science and Technology.
\end{blanko*}

%%%%%%%%%%%%%%%%%%%%%%%%%%%%%%%%%%%%%%%%%%%%%%%%%%
\section{The classical notion: \LR{} units}
%%%%%%%%%%%%%%%%%%%%%%%%%%%%%%%%%%%%%%%%%%%%%%%%%%
\label{Sec:classical}

\begin{blanko}{Semi-monoidal categories.} 
%   (B\'enabou~\cite{Benabou:CR1963}, Mac
%   Lane~\cite{MacLane:naturalAssociativity}.)  
  A {\em category with a multiplication}, or a {\em semi-monoidal
  category}, is a category $\CC$ equipped with an associative functor
  $\CC\times\CC\to\CC$, here denoted by plain juxtaposition, $(X,Y)
  \mapsto XY$.  For simplicity we assume strict associativity,
  $X(YZ)=(XY)Z$.  This is really no loss of generality: all the
  arguments in this note carry over to the case of non-strict
  associativity --- just insert associators where needed.
  
  If $X$ is an object we use the same symbol $X$ for the identity 
  arrow of $X$.
\end{blanko}

\begin{blanko}{Monoidal categories.}\label{monoidal}
  A semi-monoidal category $\CC$ is a {\em monoidal category} when it is
  furthermore equipped with a distinguished object $I$, called the 
  unit, and natural
  isomorphisms
\begin{diagram}[w=4ex,tight]
I X & \rTo^{\lambda_X}  & X 
  &  \lTo^{\rho_X}  &  X I 
\end{diagram}
obeying the following rules (cf.~Mac 
Lane~\cite{MacLane:naturalAssociativity}):
\begin{equation} 
  % Saavedra (1.3.1.1)
  % Mac Lane (5.6)
  % Kelly (4)
  \lambda_I = \rho_I
\end{equation}
\begin{equation} 
  % Saavedra (2.2.1.1 (b))
  % Mac Lane (5.3 (i))
  % Kelly (5)
  \lambda_{XY} = \lambda_X Y
\end{equation}
\begin{equation} 
  % Saavedra (2.2.1.1 (a))
  % Mac Lane (5.7 (ii))
  % Kelly (7)
 \rho_{XY} = X\rho_Y
\end{equation}
\begin{equation}\label{Kelly} 
  % Saavedra (2.2.1.1 (c))
  % Mac Lane (5.7 (i))
  % Kelly (6)
  X\lambda_Y = \rho_X Y
\end{equation}
% This object $I$ is called the unit of $\CC$.  More specifically we call
% the triple $(I,\lambda,\rho)$ an {\em LF unit}.
\end{blanko}

It was observed by Kelly~\cite{Kelly:MacLanesCoherence} that
Axiom~(\ref{Kelly}) implies the other three axioms.
We quickly run through the arguments --- they are really simple.

\begin{blanko}{Naturality.}
  Naturality of the left constraint $\lambda$ with respect to the arrow 
  $\lambda_X:IX\to X$ is expressed by this commutative diagram:
\begin{diagram}[w=6ex,h=4.5ex,tight]
IIX & \rTo^{\lambda_{IX}}  & IX  \\
\dTo<{I\lambda_X}  &    & \dTo>{\lambda_X}  \\
IX  & \rTo_{\lambda_X}  & X
\end{diagram}
Since $\lambda_X$ is invertible, we conclude
\begin{equation}\label{lambdaIX}
  \lambda_{IX} = I\lambda_X  .
\end{equation}
Similarly with the right constraint:
\begin{equation}\label{rhoXI}
  \rho_{XI} = \rho_X I   .
\end{equation}
\end{blanko}

\begin{blanko}{Fundamental observation.}\label{obs.fund}
  {\em Tensoring with $I$ from either side is an equivalence of categories:}
  $$
  I \tensor \_ \; : \; \CC \isopil \CC
  $$
  $$
  \_ \tensor I \; : \; \CC \isopil \CC
  $$

  Indeed, the left and right constraints $\lambda$ and $\rho$ 
  are invertible natural
  transformations between these two functors and the identity
  functors.  This observation does not depend on the axioms (1)--(4).
\end{blanko}

\begin{lemma}\label{Kelly-lemma}
  (Kelly~\cite{Kelly:MacLanesCoherence}.)  Axiom~(4) implies axioms
  (1), (2), and (3).
\end{lemma}
\begin{dem}
    (4) implies (2): Since tensoring with $I$ on the left is an
    equivalence, it is enough to prove $I \lambda_{XY} = I\lambda_X
    Y$.  But this equation follows from Axiom~(4) applied twice (swap
    $\lambda$ out for a $\rho$ and swap back again only on the nearest
    factor):
    $$
    I\lambda_{XY} = \rho_I XY = I \lambda_X Y .
    $$
    Similarly for $\rho$, establishing (3).
    
  (4) and (2) implies (1): Since tensoring with $I$ on the right is an
  equivalence, it is enough to prove $\lambda_I I = \rho_I I$.  But
  this equation follows from (2), (5), and (4):
  $$
  \lambda_I I = \lambda_{II} = I \lambda_I = \rho_I I .
  $$
\end{dem}

It was observed by Saavedra~\cite{Saavedra} that conversely (1), (2),
and (3) imply (4).  This will be an immediate corollary of the 
results in the next section, cf.~\ref{123->4}.

\begin{blanko}{The category of \LR{} units}
  % Saavedra (1.3.1.4)
  A triple $(I, \lambda,\rho)$ as in \ref{monoidal} is called an {\em
  \LR{} unit structure} on $(\CC,\tensor)$, or just an {\em \LR{}
  unit}.  These form the objects of a category $U(\CC)$, the {\em
  category of \LR{} units of $\CC$}, where an arrow $(I,\lambda,\rho)
  \to (I',\lambda',\rho')$ is given by an arrow $\psi: I \to I'$
  compatible with the left and right constraints in the sense that for
  every object $X$, these two triangles commute:
  \begin{equation}\label{unitmap}
  \begin{diagram}[w=6ex,h=2.5ex,tight]
  IX &&\\
  &\rdTo^{\lambda_X}  & \\
  \dTo<{\psi X}  &    & X  \\
  & \ruTo_{\lambda'_X}  &\\
  I'X  &&
  \end{diagram}
  \qquad\qquad
  \begin{diagram}[w=6ex,h=2.5ex,tight]
  && XI\\
  & \ldTo^{\rho_X}  & \\
  X && \dTo>{X \psi} \\
  & \luTo_{\rho'_X} & \\
  && X  I' .
  \end{diagram}
  \end{equation}
\end{blanko}

Note that the morphisms of \LR{} units are precisely those with respect
to which $\lambda$ and $\rho$ are also natural in $U$.  That is, one
can think of $\lambda$ (encoding all left constraints) as a natural
transformation
\vspace{15pt}
\begin{diagram}[w=3.7em]
U(\CC) \times \CC & \topglob{\tensor} \lift{-2}{\ \ \Downarrow\lambda} 
\dropglob{\text{\footnotesize{proj.}}} 
& \CC   .
\end{diagram}

\vspace{25pt}

\noindent 
The component of $\lambda$ on $((I,\lambda,\rho), X)$ is simply
$\lambda_X$.  

% \begin{MOREDETAILS}
% Note that the condition of the definition implies the 
% a priori stronger condition that for every arrow $f: X\to Y$ we have
% \begin{diagram}[w=6ex,h=4.5ex,tight]
% IX & \rTo^{\lambda^I_X}  & X  \\
% \dTo<{\phi f}  &    & \dTo>f  \\
% I'Y  & \rTo_{\lambda'_Y}  & Y
% \end{diagram}
% This is just to say that naturality in two variables follows from 
% naturality in each variable\ldots
% \end{MOREDETAILS}

\bigskip

It is a basic observation that for a semi-monoid the unit is unique if
it exists.  The following result, which goes back to 
Saavedra~\cite{Saavedra}, shows that the \LR{} unit for a monoidal
category is unique up to unique isomorphism (if it exists).

\begin{prop}\label{clasUcontract}
  The category $U(\CC)$ is contractible.
\end{prop}

\begin{dem}
  Given two \LR{} units $(I,\lambda,\rho)$ and 
  $(I',\lambda',\rho')$, define a map $\psi:I\to I'$ by
  \begin{diagram}
    I &\rTo^{\rho'_I{}^{-1}} & I I' & \rTo^{\lambda_{I'}} & I' .
  \end{diagram}
  To see that $\psi$ is compatible with the left constraint as in
  (\ref{unitmap}), consider the diagram
  \newcommand{\psimap}{%
% \setlength{\unitlength}{0.003\DiagramCellWidth}
% \multiply \unitlength by \scaleFactor
\begin{picture}(0,0)(-4,-4)
\qbezier(20,-55)(0, -55)(0, 0)
\qbezier(20,55)(0, 55)(0, 0)
\put(-17,0){\makebox(0,0)[t]{$\scriptstyle \psi X$}}
\put(20.3,-55){\vector(4,-1){0}}
\end{picture}
}

  \begin{diagram}[w=6.5ex,h=5ex,tight]
    &IX & \rTo^{\lambda_X} & X \\
    &\uTo<{\rho'_I X}>{I \lambda'_X} & & \uTo>{\lambda'_X}  \\
   \psimap &II'X & \rTo_{\lambda_{I'X}} & I'X \\
    &\dTo<{\lambda_{I'}X} && \dTo>{\lambda'_X} \\
    &I'X & \rTo_{\lambda'_X} & X
  \end{diagram}
  The bottom square is Equation~(2) composed with $\lambda'_X$.  The top square
  is naturality of $\lambda$ with respect to $\lambda'_X$.  Since all
  the arrows are invertible, we can invert the vertical arrows of the 
  top square; then the outline of the diagram is precisely the 
  left-hand compatibility triangle in (\ref{unitmap}).
  
  Compatibility with the right constraint is analogous to establish.
%   \begin{diagram}[w=3ex,h=5ex,tight]
%     &&I'X&&&\rTo^{\lambda'_X}&&&X&&&\lTo^{\rho_{X'}}&&&XI'&&\\
%     &\ruTo<{\lambda_{I'}X}&&&(2)&&&\ruTo<{\lambda'_X}&&
%                        \luTo>{\rho_X}&&&(N)&&&\luTo<{\rho_XI'}>{X\lambda_{I'}}&\\
%     II'X&&&\rTo_{\lambda_{I'X}}&&&I'X&&&&XI&&&\lTo^{\rho'_{XI}}&&&XII'\\
%     &\rdTo<{\rho'_I X}>{I \lambda'_X}&&&(N)&&&\rdTo<{\lambda'_X}&&
%                                       \ldTo>{\rho_X}&&&(2)&&&\ruTo>{X\rho'_{I}}&\\
%     &&IX&&&\rTo_{\lambda_X}&&&X&&&\lTo_{\rho_{X}}&&&XI&&
%   \end{diagram}
%   All the arrows are isomorphisms.
%   The `outer' edges of this diagram are
%   \begin{diagram}[w=6ex,h=4.5ex,tight]
%   I'X & \rTo^{\lambda'_X}  & X &  \lTo^{\rho_{X'}}&XI' \\
%   \uTo<{\psi X}  & & \uTo>X &    & \uTo>{X\psi}  \\
%   IX & \rTo_{\lambda_X}  & X &  \lTo_{\rho_{X}}&XI
%   \end{diagram}
%   which is what had to be established.
  
  Finally we check that this $\psi$ is the only morphism of \LR{} units from $I$
  to $I'$: suppose we had another morphism of \LR{} units $\gamma : I \to I'$.
  Now take $X=I$ in the left-hand compatibility
  diagram~(\ref{unitmap}).  Since all the arrows involved are
  invertible, there can by at most one compatible arrow $II \to I'I$,
  so we have $\psi I = \gamma I$.  But since tensoring with $I$ on the
  right is an equivalence, we conclude that $\psi=\gamma$.
\end{dem}

\begin{blanko}{Tensor product of \LR{} units.}\label{tensor-unit}
  There is a natural tensor on $U(\CC)$:
  the tensor product of $(I,\lambda,\rho)$ with $(I',\lambda',\rho')$
  is the object $II'$ equipped with left and right constraints
  given by these two composites:
  $$
  II'X \stackrel{I\lambda'_X}{\rTo} IX \stackrel{\lambda_X}{\rTo} X
  $$
  $$
  XII' \stackrel{\rho_X I'}{\rTo} XI' \stackrel{\rho'_X}{\rTo} X    .
  $$
  It is easy to see that this tensor product is associative.  It
  remains to check Axiom~(4) for these new left and right constraints:
  in the diagram
\begin{diagram}[w=6ex,h=4.5ex,tight]
&& XII'Y && \\
& \ldTo<{XI\lambda'_Y} && \rdTo>{\rho_X I' Y}  \\
XIY  &&&& XI'Y \\
& \rdTo_{X\lambda_Y}^{\rho_X Y} && \ldTo_{\rho'_X 
Y}^{X\lambda'_Y} & \\
&& XY &&
\end{diagram}
the inner labels show that the diagram commutes.  The passage from
inner labels to outer labels is just Axiom~(4) for
the original two units, and with the outer labels the diagram is
precisely Axiom~(4) for the new unit $II'$.

Note that the natural forgetful functor $U(\CC) \to \CC$ preserves 
the tensor product strictly.
\end{blanko}

We shall come back to this construction and show that the association
$\CC\mapsto U(\CC)$ is functorial: a strong monoidal functor
$\CC\to\DD$ induces a strong multiplicative functor $U(\CC)\to U(\DD)$
compatible with the forgetful functors (Corollary~\ref{U-functorial}).

% \begin{blanko}{Historical remark.}
%   The fact that the category of units is contractible is 
%   folklore.  Its first appearance in print is probably Saavedra's 
%   thesis~\cite{Saavedra}, Prop.~2.2.4.1.
% \end{blanko}

%%%%%%%%%%%%%%%%%%%%%%%%%%%%%%%%%%%%%%%%%%%%%%%%%%
\section{Saavedra units}
%%%%%%%%%%%%%%%%%%%%%%%%%%%%%%%%%%%%%%%%%%%%%%%%%%
\label{Sec:Saavedra}

\begin{blanko}{Cancellable objects.}
  % Saavedra, (0.1.3):  1-regular.  
  % 2-regular when the two functors are equivalences.
  An object $I$ in $\CC$ is called {\em cancellable} if the two
  functors $\CC\to\CC$
  \begin{eqnarray*}
    X & \longmapsto & IX\\ 
    X & \longmapsto & XI 
  \end{eqnarray*}
  are fully faithful.
\end{blanko}

\begin{blanko}{Idempotents.}
 A {\em idempotent} in $\CC$ is an object $I$ equipped with an
 isomorphism $\alpha : II \to I$ which is associative, i.e., the
 equation $I \alpha = \alpha I$ holds.  By a {\em pseudo-idempotent}
 we mean the same thing but without requiring $\alpha$ to be
 associative.
\end{blanko}

\begin{blanko}{Saavedra units.}\label{unit1}
  We define a {\em Saavedra unit} to be a cancellable
  pseudo-idempotent, and proceed to compare this notion with the
  notion of \LR{} unit.
\end{blanko}
The two conditions in the definition of Saavedra unit 
strengthen each other mutually:

\begin{lemma}\label{ass+eq}
  (i) A Saavedra unit $(I,\alpha)$ is in fact an idempotent,
  i.e.~$\alpha$ is automatically associative.  In other words,
  $(I,\alpha)$ is a semi-monoid.

  (ii) The two functors
  $\CC\to\CC$ given by $X\mapsto IX$ and $X\mapsto XI$ are
  in fact automatically equivalences.
\end{lemma}
This will be an easy consequence of the comparison with the
\LR{} notion of unit expressed by the next two lemmas.

\begin{lemma}\label{Saatoclas}\label{construction}
  Given a Saavedra unit $(I,\alpha)$, for each object $X$ there are
  unique arrows
\begin{diagram}[w=4ex,tight]
I X & \rTo^{\lambda_X}  & X 
  &  \lTo^{\rho_X}  &  X I 
\end{diagram}

\pagebreak

\noindent
such that 
\begin{equation}
  % Saavedra (2.2.3.3)
  \label{L}  \framebox[80pt][c]{\rule[-3pt]{0pt}{15pt} 
  $I\lambda_X \;=\; \alpha X$}
\end{equation}
\begin{equation}
  % Saavedra (2.2.3.2)
  \label{R} \framebox[80pt][c]{\rule[-3pt]{0pt}{15pt} 
  $\rho_X I \;=\; X \alpha$}
\end{equation}
The $\lambda_X$ and $\rho_X$ are isomorphisms and natural in $X$.
\end{lemma}
These two equations will be invoked throughout.

\begin{dem}
  Let $\mathbb{L}: \CC\to \CC$ denote the functor defined by
  tensoring with $I$ on the
  left.  That $\mathbb{L}$ is fully faithful means we have a bijection
  on hom sets
  $$
  \CC(IX,X) \to \CC(IIX,IX).
  $$  
  Now take $\lambda_X$ to be the 
  inverse image of $\alpha X$; it is an isomorphism since $\alpha X$
  is.  Naturality follows by considering more generally the bijection
  $$
  \Nat(\mathbb{L},\id_\CC) \to \Nat(\mathbb{L}\circ \mathbb{L}, \mathbb{L}) 
  ;
  $$
  let $\lambda$ be the inverse image of the natural
  transformation whose components are $\alpha X$.

  Similarly on the right.
\end{dem}

\begin{prop}\label{Kelly-dim1}
  The $\lambda$ and $\rho$ constructed from $\alpha$ satisfy Axiom~(\ref{Kelly}):
  $$
  \rho_X Y = X \lambda_Y .
  $$
  Hence a Saavedra unit has a canonical \LR{} unit structure.
\end{prop}
% Consequently, (1), (2), and (3) hold too.

\begin{dem}
  In the commutative square
  \begin{diagram}[w=6ex,h=4.5ex,tight]
  XIIY & \rTo^{\rho_X I Y}  & XIY  \\
  \dTo<{XI\lambda_Y}  &    & \dTo>{X\lambda_Y}  \\
  XIY  & \rTo_{\rho_X Y}  & XY
  \end{diagram}
  the left-hand arrow is equal to $X  \alpha  Y$, by $X$ tensor
  (\ref{L}),
  and the top arrow is also equal to $X  \alpha  Y$, by (\ref{R}) tensor $Y$.
  Since $X\alpha Y$ is an isomorphism, it follows that $\rho_X Y = X \lambda_Y$.
%   Precomposing with the inverse of this arrow yields the result.
\end{dem}

\begin{bevis*}{Proof of \ref{ass+eq}.}
  Re (i): set $X=Y=I$ in the Kelly equation~(4), and apply (\ref{L}) and
  (\ref{R}) once again:
  \begin{equation}\label{Kellycheck}
  I\alpha \ \stackrel{(\ref{R})}{=} \ \rho_I I \ 
  = \ I\lambda_I \ \stackrel{(\ref{L})}{=} \
  \alpha I .
  \end{equation}
  Re (ii): $\mathbb{L}$ is an equivalence, because it is
  isomorphic to the identity via $\lambda$.  (Similarly on the right.)
\end{bevis*}

Conversely, 

\begin{lemma}\label{clastoSaa}
  An \LR{} unit $(I,\lambda,\rho)$ becomes a Saavedra unit by 
  taking $\alpha \df \lambda_I = \rho_I$.
\end{lemma}
Indeed, $I$ is cancellable by \ref{obs.fund}.

\begin{blanko}{Back and forth.}
  It is clear that starting from a Saavedra unit $(I,\alpha)$ and
  constructing $\lambda$ and $\rho$ as in \ref{Saatoclas}, then
  $\lambda_I=\alpha$ and $\rho_I=\alpha$ (this follows from
  (\ref{Kellycheck})).  Conversely, starting from an \LR{} unit
  $(I,\lambda,\rho)$ and setting $\alpha\df\lambda_I = \rho_I$, then
  we have
  $$
  \alpha X = \lambda_I X \stackrel{(2)}{=} \lambda_{IX} 
  \stackrel{(5)}{=} I\lambda_X ,
  $$
  and similarly on the right, giving $X\alpha=\rho_X I$.  Hence the
  original $\lambda$ and $\rho$ satisfy the equations defining the
  left and right constraints induced from $\alpha$.  Altogether:
\end{blanko}

\begin{prop}\label{bij}
  There is a natural bijection between \LR{} unit structures (on
  $I$) and Saavedra unit structures (on $I$).
\end{prop}
We will see in a moment that this correspondence is functorial, with
the appropriate notion of morphisms of Saavedra units.  But first let
us extract a few corollaries of the back-and-forth construction.

First, since the arguments in \ref{clastoSaa} do not depend on (4),
we get the converse of Kelly's Observation~\ref{Kelly-lemma}:

\begin{cor}\label{123->4}
  % Saavedra (2.2.3.)
  Axioms (1), (2), and (3) together imply (4).
\end{cor}

As noted above, (2) and (3) are needed to ensure that the
back-and-forth construction gives back the original $\lambda$ and
$\rho$.  In absence of (2) and (3) we just get:

\begin{cor}
  Let $(\CC,\tensor)$ be a semi-monoidal category and let $I$ be an
  object equipped with natural isomorphisms $\lambda_X : IX \isopil X$
  and $\rho_X : XI \isopil X$
  such that
  $\lambda_I = \rho_I$,
  but not required to satisfy (2), (3), or (4).
  Then there exists an \LR{} unit structure on $I$ (possibly given by 
  left and right constraints different from $\lambda$ and $\rho$).
\end{cor}
Indeed, put $\alpha = \lambda_I= \rho_I$, then clearly $(I,\alpha)$ is
a Saavedra unit, and hence there are induced left and right 
constraints which satisfy (4), but there is no reason why these
new left and right constraints should coincide with 
the original  $\lambda$ and $\rho$.

\bigskip

Finally from the fact that $\lambda$ alone determines $\alpha$
(cf.~\ref{clastoSaa}), and $\alpha$ determines both $\lambda$ and
$\rho$ (by \ref{Saatoclas}), we get this:
\begin{cor}\label{redundant}
  The left constraint and the right constraint of an \LR{} unit 
  determine each other.  In other words,
  if $(I,\lambda,\rho)$ and $(I,\lambda,\rho')$ are both \LR{}
  units then $\rho=\rho'$ (and similarly with $\rho$ fixed).
\end{cor}

This corollary can also be deduced from
contractibility of the category of units (the arguments in the proof 
of Proposition~\ref{clasUcontract}): the unique morphism of units
between them is
$\lambda_I \circ \rho_I^{-1}$, which is just the identity arrow of $I$,
by (1).  Now it follows from (\ref{unitmap}) that $\rho=\rho'$.

\begin{blanko}{Economy and absoluteness of Saavedra units.}
  The \LR{} notion of unit involves a lot of structure: two whole
  families of arrows are specified, and the Corollary reveals that this
  data is somewhat redundant.  The notion of Saavedra unit is meant to
  be as economical as possible: the only structure to be specified is
  a single multiplication map $\alpha:II \to I$, a notion already 
  central to the theory of monoidal categories.

  The notion of Saavedra unit is also more fundamental than the \LR{}
  notion in that it is an absolute notion: namely, {\em the notion of
  Saavedra unit makes sense prior to any associativity constraints},
  while in contrast the axioms for an \LR{} unit only make sense
  relative to a specified associativity constraint.  In this note, for
  simplicity, the associativity is assumed strict, and in particular
  fixed, so this remark does not play any role here.  But as observed,
  all the arguments generalise to the non-strict case by inserting
  associators when needed.  Then it makes sense to fix the Saavedra
  unit and vary the associator; $\lambda$ and $\rho$ will then vary
  with the associator.  (In the non-strict setting, the associator is
  involved in Construction~\ref{construction}, and the pentagon
  equation is required to establish Axiom~(\ref{Kelly}) in
  Proposition~\ref{Kelly-dim1}).
\end{blanko}
  
\begin{blanko}{Historical remarks.}\label{error}
%   (\cite{Saavedra}, Ch.~I, 1.3.2.  Cancellable objects
%   are called {\em $1$-regular} in \cite{Saavedra} ({\em $2$-regular} meaning the
%   same thing but with equivalences instead of just fully faithful).
%   He observes that $1$-regular plus pseudo-idempotent implies 
%   $2$-regular, but does not notice that $\alpha$ is automatically 
%   associative.
  The notion of Saavedra units goes back to Saavedra's
  thesis~\cite{Saavedra} where it is mentioned in the preliminary
  chapter on monoidal categories as an alternative to the standard
  notion of unit.  He doesn't really exploit the notion, though.  He
  states the comparison result \ref{bij}, but his proof (p.34--37)
  seems to contain an error.  (First (Ch.~I,~1.3), he defines a unit
  to be a triple $(I,\lambda,\rho)$ such that $\lambda_I = \rho_I$.
  Then he imposes conditions of compatibility with the associativity
  constraint (Ch.~I,~2.2.1.1); these conditions are precisely (2),
  (3), and (4) above.  He then proves (Ch.~I,~Prop.~2.2.3) the
  converse of Kelly's result, namely that (1)+(2)+(3) imply (4).  Note
  that condition (1) can not be omitted.  The definition of Saavedra
  unit is given in Ch.~I,~1.3.2, and our Proposition~\ref{Kelly-dim1}
  is part of his Proposition~2.2.5.1.  The problem is this: after
  constructing $\lambda$ and $\rho$, he promises first to establish
  (2), (3), (4) before finally proving (1).  However to prove (4) he
  uses Prop.~2.2.3 which crucially relies on (1) (and when he comes to
  proving (1) he uses (4)).)
\end{blanko}
  
\begin{blanko}{Remark on strict units.}
  An \LR{} unit is {\em strict} if $\lambda_X$ and $\rho_X$ are
  identity arrows for all $X$.  A Saavedra unit is {\em strict} if
  $\alpha$ is the identity arrow and if the two functors `tensoring
  with $I$' are isomorphisms of categories.  To see that the latter
  strictness implies the former, note that every object $X$ is of form
  $IY$ for some $Y$.  Now $\lambda_{IY}=\alpha Y$ is an identity
  arrow, so by naturality $\lambda_X$ is an identity arrow too.
\end{blanko}

\begin{blanko}{Tensor cancellable arrows.}
  Let $I$ and $J$ be cancellable objects.  An arrow $\psi: I \to J$ is
  called {\em left tensor cancellable} if the induced map on hom sets
  \begin{eqnarray}\label{Homphi}
    \CC(X,Y) & \longrightarrow & \CC(IX,JY)  \\
    f & \longmapsto & \psi f \notag
  \end{eqnarray}
  is a bijection.  Right tensor cancellable is defined in the obvious
  way; an arrow is {\em tensor cancellable} if it is both left and
  right cancellable.
\end{blanko}

\begin{blanko}{Morphisms of Saavedra units.}
  % Saavedra 1.3.2.2 just requires them to be isomorphisms, without
  % trying to explain why\ldots
  Let $(I,\alpha)$ and $(J,\beta)$ be Saavedra units in $\CC$.  
  In partcular they are semi-monoids, by \ref{ass+eq}~(i).  A {\em
  morphism of Saavedra units} is a tensor cancellable semi-monoid
  homomorphism $\psi:I\to J$.  Being a semi-monoid homomorphism means
  that this diagram commutes:
  \begin{equation}\label{semi-monoid-homo}
    \begin{diagram}[w=6ex,h=4.5ex,tight]
    II & \rTo^{\psi\psi}  & JJ  \\
    \dTo<\alpha  &    & \dTo>\beta  \\
    I  & \rTo_\psi  & J .
    \end{diagram}
  \end{equation}
  This defines the category of Saavedra units in $\CC$.
\end{blanko}

\begin{BM}
  By factoring (\ref{Homphi}) in two ways like this:
  \begin{diagram}[w=10ex,h=6ex,tight]
  \CC(X,Y) & \rTo^{J \CC }  & \CC(JX,JY)  \\
  \dTo<{I\CC}  &    & \dTo>{\text{precomp. } \psi X}  \\
  \CC(IX,IY)  & \rTo_{\text{postcomp. } \psi Y}  & \CC(IX,JY)
  \end{diagram}
  it follows that any monomorphism or epimorphism between
  cancellable objects is tensor cancellable.  In particular, isomorphisms 
  are tensor cancellable.
  
  It follows from the next Proposition that a morphism of Saavedra
  units is automatically an isomorphism, since morphisms of \LR{}
  units are so:
\end{BM}

\begin{prop}\label{unit-morphisms}
  The category of Saavedra units is canonically isomorphic to the 
  category of \LR{} units, hence contractible.
\end{prop}
% This is Saavedra 2.2.5.1, but with a wrong proof.

% Saavedra considers also units which are not required to be
% compatible with associativity.  In that case the category of
% \LR{} units is not necessarily connected, and the functor from
% \LR{} units to Saavedra units is not necessarily an isomorphism.
% Saavedra claims it is fully faithful however...

\begin{dem}
  We have already established a bijection on the level of objects.
  
  Given an arrow in the category of \LR{} units, $\psi: I \to J$,
  then in particular it is an isomorphism and hence tensor cancellable.  It
  remains to check that it is a semi-monoid homomorphism $(I,\alpha)
  \to (J,\beta)$.  But this is easy: in the diagram
  
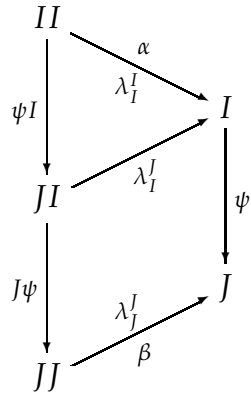
\begin{diagram}[w=6ex,h=3ex,tight]
  II &&  \\
    & \rdTo_{\lambda^I_I}^\alpha   &   \\
  \dTo<{\psi I}  &   & I \\
  & \ruTo_{\lambda^J_I} \\
  JI && \dTo>\psi \\
  \\
  \dTo<{J\psi} && J \\
  & \ruTo_\beta^{\lambda^J_J} \\
  JJ
  \end{diagram}
  the triangle is compatibility~(\ref{unitmap}) with the left
  constraint for the object $I$.  The square is naturality of
  $\lambda^J$ with respect to $\psi$.  The outer square is the
  semi-monoid homomorphism condition (\ref{semi-monoid-homo}).
  
  \bigskip
  
  Conversely, suppose we have a tensor cancellable semi-monoid homomorphism
  $\psi:I \to J$.  Construct $\lambda^I_X : IX \to X$ and $\lambda^J_X : 
  JX \to X$ as in \ref{construction}, and check that the composite
  $$
  \theta: \qquad IX \stackrel{\psi X}{\rTo} JX \stackrel{\lambda^J_X}{\rTo} X
  $$
  is equal to $\lambda^I_X$: In the diagram
    
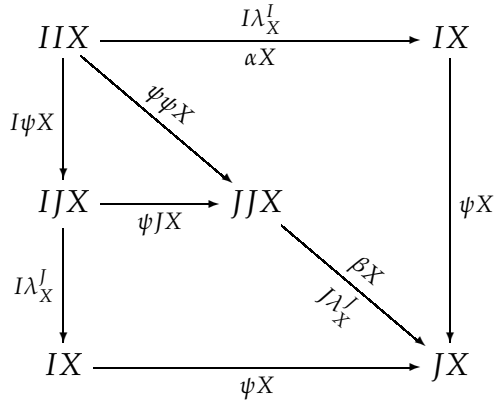
\begin{diagram}[w=6.5ex,h=5.5ex,tight,hug]
  IIX  && \rTo^{I\lambda^I_X}_{\alpha X} && IX\\
   \dTo<{I\psi X} &\rdTo^{\psi \psi X} & && \\
   IJX &\rTo_{\psi JX} & JJX && \dTo>{\psi X} \\
   \dTo<{I\lambda^J_X} &&&\rdTo^{\beta X}_{J\lambda^J_X}&\\
  IX && \rTo_{\psi X} && JX
  \end{diagram}
  the upper right-hand cell is the semi-monoid homomorphism condition
  (\ref{semi-monoid-homo}) (tensored with $X$ on the right), and the
  other cells are trivially commutative.  Now the left-and-bottom
  composite is the tensor product $\psi \theta$ while the
  top-and-right composite is $\psi \lambda^I_X$.  These two coincide
  since the outer square commutes.  But since $\psi$ is (left) tensor
  cancellable, we conclude that $\theta=\lambda^I_X$ as wanted.  A
  similar argument shows that $\psi$ is compatible with the right
  constraint.
\end{dem}

Observe that in the passage from morphism of \LR{} units to
morphism of Saavedra units we only used compatibility with the left
constraint (and could equally well have used only $\rho$), 
and in fact we only used the compatibility with respect to
the object $X=I$.  (It would also have been enough to test
the compatibility with respect to the object $X=J$.)
In particular:

\begin{cor}\label{compLR}
   Given \LR{} units $(I,\lambda,\rho)$ and $(I',\lambda',\rho')$,
   an arrow $I \to I'$ is compatible with the left constraints if and
   only if it is compatible with the right constraints.
\end{cor}

\begin{blanko}{Tensor products of Saavedra units.}
  The isomorphism of Proposition~\ref{unit-morphisms} endows the
  category of Saavedra units with a tensor product, which is
  surprising if you only look at the definition of Saavedra unit,
  where the only structure is a semi-monoid structure --- usually in
  order to define tensor products of semi-monoids you need a symmetry
  on the underlying semi-monoidal category.

  Tracing through the correspondences, the tensor product of two 
  Saavedra units $(I,\alpha)$ and $(J,\beta)$ is given by
  $(IJ,\gamma)$ where $\gamma$ is the composite
  $$
  IJIJ \stackrel{IJ\lambda^I_J}{\rTo}  IJJ \stackrel{\rho^J_I J}{\rTo} 
  IJ .
  $$
\end{blanko}

%%%%%%%%%%%%%%%%%%%%%%%%%%%%%%%%%%%%%%%%%%%%%%%%%%
\section{Strong functoriality}
%%%%%%%%%%%%%%%%%%%%%%%%%%%%%%%%%%%%%%%%%%%%%%%%%%
\label{Sec:functoriality}

In this section we prove that the category of monoidal categories and
strong monoidal functors in the \LR{}-unit sense is isomorphic
to the category of monoidal categories and strong monoidal functors in
the Saavedra-unit sense.

\begin{blanko}{Strong multiplicative functors.}
  A {\em strong multiplicative functor} is a functor
  \begin{eqnarray*}
      \CC & \longrightarrow & \DD  \\
      X & \longmapsto & \OV{X}
  \end{eqnarray*}
  equipped with a natural family of isomorphisms in $\DD$:
  $$
  \phi_2\df\phi_{X,Y} : \OV{X} \OV Y \to \OV{XY}
  $$
  such that this square commutes:
  \begin{equation}\label{phi2-ass}
  \begin{diagram}[w=6ex,h=4.5ex,tight]
  \OV{X} \OV Y \OV Z & \rTo^{}  & \OV{XY}\OV Z  \\
  \dTo  &    & \dTo  \\
  \OV{X} \OV{YZ}  & \rTo  & \OV{XYZ}
  \end{diagram}
  \end{equation}
%   The identity functor of a semi-monoidal category canonically 
%   becomes semi-monoidal by taking the identity comparison maps,
%   and the composite of two semi-monoidal functors canonically
%   acquires the structure of a semi-monoidal functor.
%   
%   We will mostly be concerned with the case the comparison arrows are
%   required to by invertible, in which case we talk about strong
%   multiplicative functors.  If the comparison arrows are identity
%   arrows we say the multiplicative functor is strict.
\end{blanko}

%   \begin{MOREDETAILS}
%     Given $\CC\to\DD\to\EE$ denoted by bar and tilde respectively,
%     $X \mapsto \OV{X} \mapsto \tilde{\OV{X}}$, and if they are 
%     semi-monoidal $\phi:\OV{X} \OV Y \to \OV{XY}$ and $\psi:\tilde P \tilde Q 
%     \to \wtil{PQ}$, then the composite has this map as structure map
%     $$
%     \tilde{\OV{X}} \tilde {\OV Y} \stackrel{\psi}{\rTo}
%     \wtil{\OV{X} \OV{Y}} \stackrel{\tilde \phi}{\rTo}
%     \wtil{\OV{XY}}
%     $$
%   \end{MOREDETAILS}
% 

\begin{blanko}{Strong monoidal functors in the \LR{}-unit sense.}
  \label{strong-clas}
  A {\em strong monoidal functor} is a strong multiplicative functor
  $(\CC,\tensor,I)\to(\DD,\tensor,J)$, $X\mapsto \OV{X}$ together with
  an isomorphism $\phi_0: J \to \stregI$ satisfying these two
  conditions of compatibility with the left and right constraints:
  \begin{equation}\label{strong-clas-comp}
  \begin{diagram}[w=6ex,h=4.5ex,scriptlabels,tight]
  J \OV{X}       &\rTo^{\lambda^J_{\OVscriptscript{X}}}    & \OV{X}    \\
  \dTo<{\phi_0 \OVscript{X}}    &      & \uTo_{\OVscript{\lambda^I_X}}    \\
  \stregI \OV{X}&   \rTo_{\phi_2} &  \OV{IX}
  \end{diagram}
  \hspace{5em}
  \begin{diagram}[w=6ex,h=4.5ex,scriptlabels,tight]
  \OV{X}    & \lTo^{\rho^J_{\OVscriptscript{X}}}    & \OV{X}  J   \\
  \uTo<{\OVscript{\rho^I_X}}    &        & \dTo>{\OVscript{X}\phi_0} \\
   \OV{XI}    & \lTo_{\phi_2}  & \OV{X} \stregI
  \end{diagram}
  \end{equation}
  In Section~\ref{Sec:lax} we shall consider also lax monoidal 
  functors, where $\phi_2$ and $\phi_0$ are not required to be 
  invertible.
\end{blanko}

\begin{BM}
  Saavedra~\cite{Saavedra} considers only strong monoidal functors
  (and only in the \LR{}-unit sense): he requires $\phi_2$ to be an
  isomorphism, and claims that then $\phi_0$ is automatically an
  isomorphism too (\cite{Saavedra}, Ch.~I,~4.2.3.).  This claim is
  false.  It is true that $\phi_0 \OV{X}$ is always an isomorphism, as
  seen in the diagram.  But if $\OV{X}$ is not a cancellable object in
  $\DD$, then this does not imply that $\phi_0$ is an isomorphism.  As
  a concrete counter example, consider the lax monoidal functor
  $\Vect\to\Vect$ sending every vector space to the zero vector space.
  This is a strong multiplicative functor (in fact strict), but the
  unit comparison $\phi_0 : \ground \to \ov{\ground}=\{0\}$ is
  obviously not an isomorphism.
  In the following paragraph, Saavedra states that a monoidal functor
  with a compatibility with given units might not have a compatibility
  with other units.  This is correct as stated there, because at that
  point he does not assume Axiom~(\ref{Kelly}) to hold.  In that case the
  category of units is not necessarily connected, and compatibility
  is only guaranteed for units connected to the given compatible unit.
  When Axiom~(\ref{Kelly}) is assumed to hold, as throughout the present 
  note, compatibility with one unit implies compatibility with any
  other unit.  This observation, and some related ones, will be made
  in the setting of Saavedra units where they are trivial to 
  establish (cf.~\ref{comp-all}).
\end{BM}

\begin{blanko}{Strong monoidal functors in the Saavedra-unit sense.}
  \label{strong-Saa}
  In the viewpoint of Saavedra units, a {\em strong monoidal functor}
  is defined to be a strong multiplicative functor
  $(\CC,\tensor,(I,\alpha))\to(\DD,\tensor,(J,\beta))$ together with a
  semi-monoid isomorphism $\phi_0 : J \to \stregI$.
  %   satisfying
  %   \begin{diagram}[w=6ex,h=4.5ex,tight]
  %   J & \rTo^{\phi_0}  & \stregI  \\
  %   \uTo<\beta  &    & \uTo>{\OVscript{\alpha}}  \\
  %   J J  & \rTo_{\phi_0  \phi_0}  & \stregI \stregI
  %   \end{diagram}
  %   (Here $\OV{\alpha}$ is really a composite  $\stregI  \stregI
  %   \to \OV{I  I} \to \stregI$ involving the multiplicative
  %   compatibility.)
  % 
  %   \bigskip
  % 
  %   So it is a multiplicative functor together with a semi-monoid
  %   isomorphism $\phi_0: J \isopil \stregI$.
\end{blanko}
  Recall that $\stregI$ is a semi-monoid via the arrow
  \begin{equation}\label{streg-alpha}
  \stregI \stregI \stackrel{\phi_2}{\rTo} \OV{II} 
  \stackrel{\OVscript{\alpha}}{\rTo} \stregI ,
  \end{equation}
  so the semi-monoid homomorphism condition is this:
  \begin{equation}\label{strong-Saa-semi}
%   \begin{diagram}[w=6ex,h=4.5ex,tight]
%   J & \rTo^{\phi_0}  & \stregI  \\
%   \uTo<{\beta}  &    & \uTo>{\OVscript{\alpha}}  \\
%   JJ  &   & \OV{II} \\
%   & \rdTo_{\phi_0 \phi_0} & \uTo>{\phi_2} \\
%   && \stregI \stregI
%   \end{diagram}
  \begin{diagram}[w=6ex,h=3.5ex,tight]
  JJ & \rTo^{\phi_0 \phi_0}  & \stregI \stregI  \\
    &    & \dTo>{\phi_2}  \\
  \dTo<{\beta}  &   & \OV{II} \\
  && \dTo>{\OVscript{\alpha}} \\
  J & \rTo_{\phi_0} & \stregI
  \end{diagram}
  \end{equation}  

\begin{prop}\label{compatibilities}
  The two compatibility conditions are equivalent.  In other words,
  for a strong multiplicative functor, a compatibility with an \LR{}
  unit canonically induces a compatibility with the corresponding
  Saavedra unit, and vice versa.
\end{prop}
For later reference, we split the statement into two lemmas,
Lemma~\ref{funct:clas->Saa} and Lemma~\ref{funct:Saa->clas}.
But first a couple of remarks, the first one being rather 
important:

\begin{prop}\label{phi-unitmorph}
  $\phi_0 : J \to \stregI$ is a Saavedra-unit compatibility if and only if
  $(\stregI,\OV{\alpha})$ is a Saavedra unit and $\phi_0$ is a
  morphism of Saavedra units.
\end{prop}
Here and below, by abuse of notation, the symbol $\OV{\alpha}$ refers
to the composite~(\ref{streg-alpha}).

\begin{dem}
  Suppose $\phi_0$ is a unit compatibility, cf.~\ref{strong-Saa}.  The
  fact that $\stregI$ is isomorphic to $J$ implies that it is itself
  cancellable, and hence $(\stregI,\OV{\alpha})$ is a Saavedra unit.
  The fact that $\phi_0$ is an invertible semi-monoid homomorphism is
  just to say that it is a morphism of units.  Conversely, if
  $(\stregI,\OV{\alpha})$ is a Saavedra unit and $\phi_0$ is a
  morphism of units, then in particular it is an invertible
  semi-monoid homomorphism.
\end{dem}

Note that this proposition can not be stated in terms of \LR{} units,
since a priori the image of an \LR{} unit does not even have enough 
structure to make sense of the question whether it satisfies the 
unit axioms.  (The left and right constraints of the image of the unit
are only defined for objects in the image of the functor.)

Combining this proposition with contractibility of the category of units
we get these three immediate corollaries:

\begin{cor}\label{unit-comp-unique}
  A unit compatibility on a strong multiplicative functor is unique if it 
  exists.
\end{cor}
This can be seen as a generalisation of the uniqueness
part of Proposition~\ref{clasUcontract}.  Namely, to give a morphism
of units $I \to I'$ in $\CC$ is the same thing as providing a unit
com\-patibility on the identity functor on $\CC$, with $I$ chosen as
unit in the domain and $I'$ in the codomain.

\begin{cor}\label{comp-all}
  Let $\CC$ and $\DD$ be monoidal categories and let 
  $\CC\to\DD$ be a strong multiplicative functor.  Given a unit compatibility
  with respect to chosen units $I$ in $\CC$ and $J$ in $\DD$, then 
  there are canonical compatibilities with any other choices of units 
  in $\CC$ and $\DD$.
\end{cor}
Finally, in view of Corollary~\ref{unit-comp-unique} we can consider
unit compatibility to be a property, not a structure, and restate
Proposition~\ref{phi-unitmorph} as:

\begin{cor}\label{image-of-unit}
  A multiplicative functor is monoidal if and only if the image of a 
  unit is again a unit.
\end{cor}
Again, this statement does not even make sense for \LR{} units.

\bigskip

Now for the lemmas that make up Proposition~\ref{compatibilities}:

\begin{lemma}\label{funct:clas->Saa}
  If $\CC\to\DD$, $X\mapsto \OV{X}$ is a multiplicative functor, and
  $\phi_0: J \to \stregI$ is an \LR{}-unit compatibility, then it
  is also a Saavedra-unit compatibility.
\end{lemma}

\begin{dem}
%   Suppose $(F,\phi_0)$ is compatible with $\lambda$ and $\rho$.
  In fact it is enough to have the 
  compatibility~(\ref{strong-clas-comp}) with $\lambda$ for the 
  object $X=I$. We then have a commutative diagram
  \begin{diagram}[w=6ex,h=4.5ex,scriptlabels,tight]
  J J       &\rTo_{\lambda^J_{J}}^{\beta} & J    \\
  \dTo<{J \phi_0}    &      & \dTo_{\phi_0}    \\
  J \stregI       &\rTo^{\lambda^J_{\stregIscriptscript}}    & \stregI    \\
  \dTo<{\phi_0  \stregIscript}    &      & 
  \uTo<{\OVscript{\lambda^I_I}}>{\OVscript{\alpha}} \\
  \stregI \stregI&   \rTo_{\phi_2} &  \OV{I I}
  \end{diagram}
  The top square is naturality of $\lambda^J$ with respect to
  $\phi_0$, and the bottom square is the left compatibility 
  (\ref{strong-clas-comp}).  The
  outer square is precisely the compatibility 
  diagram~(\ref{strong-Saa-semi}) for $\phi_0$
  with respect to the Saavedra units $\beta=\lambda^J_J$ and
  $\alpha=\lambda^I_I$.
\end{dem}

\begin{lemma}\label{funct:Saa->clas}
  If $\phi_0: J \to \stregI$ is a Saavedra-unit compatibility then it is
  also compatible with $\lambda$ and $\rho$.
\end{lemma}

\begin{dem}
  In the diagram
  \begin{equation}\label{bigdiagram}
\begin{diagram}[w=6ex,h=4.5ex,tight]
  JJ\OV{X} &&&\rTo_{\beta \OVscript{X}}^{J \lambda^J_{\OVscriptscript{X}}}
  &&& J \OV{X}\\
  \dTo<{J\phi_0 \OVscript{X}}&&&&&&\dTo>{\phi_0\OVscript{X}} \\
J \stregI \OV{X} & \rTo^{\phi_0 \stregIscript \OVscript{X}}  & \stregI \stregI \OV{X} 
& \rTo^{\phi_2 \OVscript{X}} & \OV{II} \OV{X} &\rTo^{\OVscript{\alpha} 
\OVscript{X}} & \stregI \OV{X}\\
\dTo<{J \phi_2}&    & \dTo<{\stregIscript \phi_2} && \dTo>{\phi_2} && \dTo>{\phi_2}  \\
J \OV{IX}  & \rTo_{\phi_0 \OVscript{IX}}  & \stregI \OV{IX} & \rTo_{\phi_2} 
& \OV{IIX} & \rTo^{\OVscript{\alpha X}}_{\OVscript{\lambda_I X}} & \OV{IX} \\
\dTo<{J \OVscript{\lambda^I_X}}&&\dTo<{\stregIscript \OVscript{\lambda^I_X}} && 
\dTo>{\OVscript{I\lambda^I_X}} && \dTo>{\OVscript{\lambda^I_X}} \\
J \OV{X} & \rTo_{\phi_0 \OVscript{X}} & \stregI \OV{X} & \rTo_{\phi_2} & \OV{IX} 
& \rTo_{\OVscript{\lambda^I_X}} & \OV{X}
\end{diagram}
\end{equation}
the top square is precisely the compatibility~(\ref{strong-Saa-semi})
between $\alpha$ and $\beta$ (tensored with $X$); the other squares
are commutative for trivial reasons.  Now since the bottom composite
and the right-hand composite in the big square are equal and are
monomorphisms, we can cancel them away and conclude that the left-hand
composite is equal to the top arrow as we wanted.

The compatibility with $\rho$ is established analogously.
\end{dem}

In the proof of Lemma~\ref{funct:clas->Saa} we only used
compatibility with $\lambda_I$, so as a corollary we get the following
result, which does not refer to the notion of Saavedra unit.

\begin{cor}\label{functor-redundance}
  A strong monoidal functor (in the classical sense) is compatible
  with left constraints if and only if it is compatible with right
  constraints, and this can be measured on $I$ alone.
\end{cor}

Joining Proposition~\ref{compatibilities} and Proposition~\ref{bij}, 
and with the monoidal natural transformations as 
$2$-cells (cf.~just below), we get

\begin{prop}\label{2-iso}
  There is an isomorphism 
%   of $2$-categories 
%   $\kat{MonCat} \simeq
%   \kat{RU-MonCat}$.
  between the $2$-category of monoidal categories, strong monoidal
  functors (in the \LR{} sense), and monoidal natural
  transformations, and the $2$-category of monoidal categories, strong
  monoidal functors in the Saavedra-unit sense, and monoidal natural
  transformations.
\end{prop}

\begin{blanko}{Monoidal natural transformations.}
  A natural transformation $u:(F,\phi_2)\Rightarrow (G,\gamma_2)$ between
  strong multiplicative functors $F:\CC\to\DD$, $X\mapsto \OV{X}$ and 
  $G:\CC\to\DD$, $X\mapsto \wtil X$ is called  {\em multiplicative}
  if for every pair of objects $X,Y$ in $\CC$ this diagram commutes:
    \begin{diagram}[w=6ex,h=4.5ex,tight]
    \OV{X} \OV{Y} & \rTo^{\phi_2} & \OV{XY}   \\
    \dTo<{u_X u_Y}  &    & \dTo>{u_{XY}}  \\
   \wtil X \wtil Y & \rTo_{\gamma_2} & \wtil{XY}   .
    \end{diagram}

    A multiplicative natural transformation is {\em monoidal} if the 
    following unit 
  condition is satisfied:
  \begin{equation}\label{nt-unit-cond}
    \begin{diagram}[w=3ex,h=3.5ex,tight]
      \stregI  && \rTo^{u_I}    && \wtil I    \\
    &\luTo_{\phi_0}    &      & \ruTo_{\gamma_0}  &  \\
    &    & J   & &
    \end{diagram}
  \end{equation}
  (It is the same condition for \LR{} units and Saavedra units.)
\end{blanko}

There is one remark to make about monoidal natural transformations:
\begin{lemma}
  Condition (\ref{nt-unit-cond}) holds automatically if 
  $u_I$ is an isomorphism (or just tensor cancellable).
\end{lemma}

\begin{dem}
  By \ref{phi-unitmorph}, $(\stregI,\OV{\alpha})$ and $(\wtil I,\wtil\alpha)$
  are Saavedra units, and $\phi_0$ and $\gamma_0$ are morphisms 
  of units.  Now,
  naturality of $u$ means that we have this commutative square
\begin{diagram}[w=6ex,h=4.5ex,tight]
\OV{II} & \rTo^{u_{II}}  & \wtil{II}  \\
\dTo<{\OVscript{\alpha}}  &    & \dTo>{\wtil\alpha}  \\
\stregI  & \rTo_{u_I}  & \wtil I
\end{diagram}
Hence if just $u_I$ is tensor cancellable then it is a morphism of units, and
hence the triangle~(\ref{nt-unit-cond}) commutes --- every diagram
of morphisms of units commutes because the category of units is contractible.
\end{dem}

% \begin{blanko}{Historical remarks.}
%   Saavedra does not consider monoidal functors in the viewpoint of
%   Saavedra units.  He gives a proof of the previous lemma using
%   only \LR{} units, but it is quite complicated.
% \end{blanko}

\begin{blanko}{Unit compatibility in terms of multiplicativity for units.}
  The following result is a variation of the corollaries
  \ref{unit-comp-unique} and \ref{comp-all}, but stated in a global
  manner involving also multiplicativity.  It will play an important
  role in Section~\ref{Sec:gentle}.
\end{blanko}

\begin{prop}\label{extend}
  Let $\CC\to\DD$ be a multiplicative functor admitting a unit
  compatibility.  The totality of all units compatibilities (one for
  each choice of unit in $\CC$ and in $\DD$) amounts precisely to a
  lift to a strong multiplicative functor $U(\CC)\to U(\DD)$.  I.e., a
  commutative diagram of strong multiplicative functors
  \begin{equation}\label{UCUD}
  \begin{diagram}[w=6ex,h=4.5ex,tight]
  U(\CC) & \rTo  & U(\DD)  \\
  \dTo  &    & \dTo  \\
  \CC  & \rTo  & \DD
  \end{diagram}
  \end{equation}
\end{prop}

\begin{dem}
  Given a lift, then for any $I$ in $U(\CC)$ and $J$ in $U(\DD)$,
  there is a connecting isomorphism $J \isopil \stregI$, since 
  $U(\DD)$ is contractible.  Hence by \ref{phi-unitmorph} we have a 
  unit compatibility.  
  
  The converse implication is subtler.  We already noticed
  (\ref{phi-unitmorph}) that the existence of a unit compatibility
  implies that the image of any unit is again a unit.  It is also easy
  to see that if $\psi$ is a morphism of units then $\OV{\psi}$ is a
  morphism of units --- this is just a consequence of the naturality
  of $\phi_0$ with respect to $\OV{\psi}$.  Composition is obviously
  preserved; hence we have a functor $U(\CC)\to U(\DD)$, and it is
  clear that diagram~(\ref{UCUD}) commutes as a diagram of functors.
  
  Furthermore, since $U(\DD)$ is contractible there exists a unique
  multiplicative structure on the functor $U(\CC)\to U(\DD)$.  Indeed,
  its components $\stregI \OV{I}{}' \to \OV{II'}$ are the unique
  comparison arrows that exist in $U(\DD)$, and axiom~(\ref{phi2-ass})
  is automatically satisfied since all diagrams commute in a
  contractible category.
  
  It remains to check that this map extends the $\phi_2$ of the
  original functor, i.e., that the diagram commutes as strong
  multiplicative functors, not just as functors.  Verifying this
  amounts to checking that $\phi_2 : \stregI \OV{I}{}' \to \OV{II'}$
  is a morphism of units, i.e.~a semi-monoid homomorphism.  Doing this
  in the Saavedra-unit setting is a bit cumbersome since the
  definition of the tensor products in the categories of units involve
  the left and right constraints anyway.  But it is not difficult to
  check that $\phi_2 : \stregI \OV{I}{}' \to \OV{II'}$ is a morphism
  of \LR{} units.  We need to check that this diagram commutes
  for all objects $Y$ in $\DD$:
  \begin{equation}\label{phiY}
    \begin{diagram}[w=6ex,h=2.5ex,tight]
  \stregI \OV{I}' Y &&\\
  &\rdTo^{\lambda^{\stregIscriptscript \OVscriptscript I'}_Y}  & \\
  \dTo<{\phi_2 Y}  &    & Y  \\
  & \ruTo_{\lambda^{\OVscriptscript {II}'}_Y}  &\\
  \OV{II'}Y  &&
  \end{diagram}
  \end{equation}
  In view of Corollary~\ref{compLR} we do not have to check also
  the right-hand diagram, although of course this could be done 
  analogously.   By the remarks preceding that corollary, it is 
  furthermore enough to check the diagram in the case $Y= \OV{II'}$.
  Since there are already so many $I$s involved, it is practical
  in the following
  argument to set $\OV{X}\df \OV{II'}$.
  
  Note that for the unit $(\stregI,\OV{\alpha})$, the corresponding 
  left constraint $\lambda^{\stregIscript}_{\OVscript{X}}$
  is given by the  equation
  \begin{equation}\label{leftbar}
    \begin{diagram}[w=6ex,h=2.5ex,tight]
  \stregI \OV{X} &&\\
  &\rdTo^{\lambda^{\stregIscriptscript}_{\OVscript{X}}}  & \\
  \dTo<{\phi_2}  &    & \OV{X}  \\
  & \ruTo_{\OVscript{\lambda^I_X}}  &\\
  \OV{IX}  &&
  \end{diagram}
  \end{equation}
  Now the wanted equation~(\ref{phiY}) is the outline of this 
  big diagram:
  
  
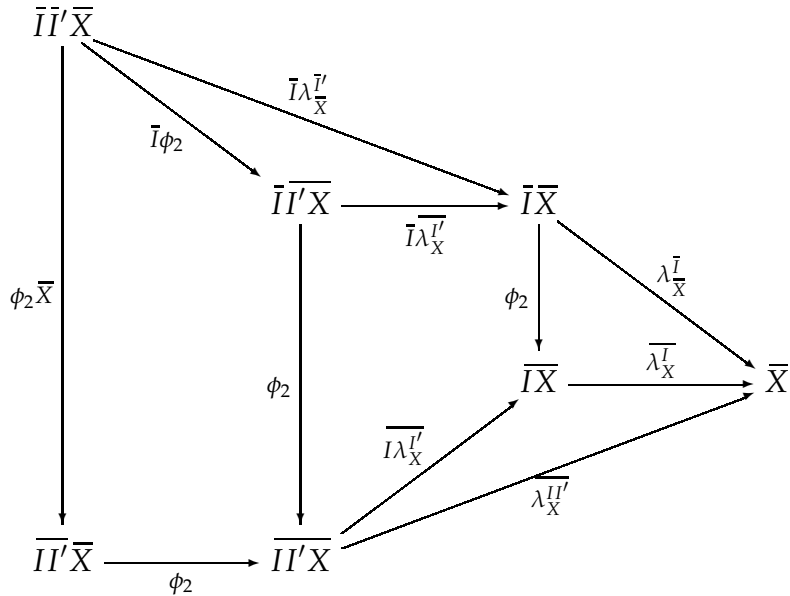
\begin{diagram}[w=8ex,h=6ex,tight]
  \stregI \OV{I}{}' \OV{X} &&&&&& \\
  & \rdTo_{\stregIscript \phi_2}  \rdTo(4,2)^{\stregIscript 
  \lambda^{\stregIscriptscript{}'}_{\OVscriptscript{X}}} &&&&&  \\
  && \stregI \OV{I'X}  & \rTo_{\stregIscript 
  \OVscript{\lambda^{I'}_X}}  & \stregI \OV{X} && \\
  \dTo<{\phi_2 \OVscript{X}} &&  && \dTo<{\phi_2} & 
  \rdTo^{\lambda^{\stregIscriptscript}_{\OVscriptscript{X}}} & \\
  &&\dTo<{\phi_2}&& \OV{IX} & \rTo^{\OVscript{\lambda^I_X}} & \OV{X} \\
  &&&\ruTo^{\OVscript{I \lambda^{I'}_X}} && 
  \ruTo(4,2)_{ \ \OVscript{\lambda^{II'}_X}} & \\
  \OV{II'}\OV{X}&\rTo_{\phi_2}&\OV{II'X} &&&&  
  \end{diagram}
  The left-hand square is associativity of $\phi_2$.  The middle 
  square is naturality of $\phi_2$.  The three triangles are 
  the definition of the left constraints as in (\ref{leftbar}).
  Now the top composite is precisely the left constraint of the 
  tensor product $\stregI\stregI{}'$ (by definition~\ref{tensor-unit}),
  and the bottom composite is
  the image of the left constraint on $II'$.
\end{dem}

\begin{cor}\label{U-functorial}
  The association $\CC\to U(\CC)$ defines a functor from the category
  of monoidal categories and strong monoidal functors to the category 
  of categories with multiplication and strong multiplicative functors.
\end{cor}

%%%%%%%%%%%%%%%%%%%%%%%%%%%%%%%%%%%%%%%%%%%%%%%%%%
\section{Monoids}
%%%%%%%%%%%%%%%%%%%%%%%%%%%%%%%%%%%%%%%%%%%%%%%%%%
\label{Sec:monoids}

In this short section we describe what a monoid is in the Saavedra-unit
setting.  The proof is postponed to the end of the paper where it will
be a special case of the treatment of lax monoidal functors.
The Saavedra-unit notion of monoid motivates the fancy description of
Saavedra units in Section~\ref{Sec:gentle}.

\begin{blanko}{Semi-monoids.}
  A {\em semi-monoid} in a semi-monoidal category $\CC$ is an object $S$
  equipped with an associative  multiplication map $\mu: SS \to S$.
  A {\em semi-monoid homomorphism} from $(S,\mu)$ to $(S',\mu')$ 
  is an arrow $S\to S'$ such that this diagram commutes:
  \begin{diagram}[w=6ex,h=4.5ex,tight]
  SS & \rTo  & S'S'  \\
  \dTo<{\mu}  &    & \dTo>{\mu'}  \\
  S  & \rTo  &S' . 
  \end{diagram}
\end{blanko}

\bigskip

Henceforth when we refer to a monoidal category, we understand it to
have both \LR{} and Saavedra unit structure, corresponding to each 
other as described in 
Section~\ref{Sec:Saavedra}.  So a monoidal category is the data of
$(\CC,\tensor, (I,\alpha), \lambda,\rho)$.

\begin{blanko}{Classical monoids.}
  A {\em monoid} in the classical sense is a semi-monoid $\mu: MM \to
  M$ equipped with an arrow $\eta: I \to M$ such that these two
  triangles commute:
  \begin{equation}\label{monoidunitcond}
    \begin{diagram}[w=6ex,h=4.5ex,tight]
    IM & \rTo^{\eta M}  & MM & \lTo^{M\eta} & MI  \\
      & \rdTo_{\lambda_M}   &\dTo>\mu& \ldTo_{\rho_M}  \\
      &   & M
    \end{diagram}
  \end{equation}
  The arrow $\eta:I\to M$ with these properties is unique,
  if it exists.
\end{blanko}

% \begin{blanko}{One step backwards: monoids.}
%   The alternative axioms for monoidal categories can also be applied to
%   monoids in $\Set$: a monoid can be defined as a semi-monoid $C$
%   equipped with an object $I$ such that $II=I$ and such the maps $C
%   \to C$ defined by multiplication by $I$ (from either side) are
%   bijections.
% \end{blanko}

\begin{blanko}{Gentle maps.}
  Let $(M,\mu)$ be a semi-monoid, and let $U$ be any object.
  We call an arrow $U \to M$ {\em gentle}
  if the two composites
  $$
  UM \to MM \stackrel{\mu}{\to} M
  $$
  $$
  MU \to MM \stackrel{\mu}{\to} M
  $$
  are monomorphisms.
\end{blanko}

\begin{blanko}{Saavedra monoids.}\label{Saavedra-monoid}
  A {\em Saavedra monoid} is a semi-monoid $\mu : MM \to M$
  equipped with a gentle
  semi-monoid homomorphism
$\eta: (I,\alpha) \to (M,\mu)$  (i.e.~a semi-monoid homomorphism 
that is gentle as an arrow in $\CC$.)

As a consequence of the following proposition, in fact the two
monomorphisms are automatically isomorphisms.
\end{blanko}

\begin{prop}\label{monoid-prop}
  There is an isomorphism between the category of 
  classical monoids in $\CC$ and the category of Saavedra monoids in 
  $\CC$.
\end{prop}

The arrows in these categories are described in the following
paragraph.  The Proposition can be proved directly without difficulty,
but we postpone the proof --- it will be a special case of the
treatment of lax monoidal functors, cf.~page~\pageref{monoid-proof}.

\begin{blanko}{Monoid homomorphisms.}
  A {\em monoid homomorphism} between two monoids in $\CC$ is just a 
  semi-monoid homomorphism $\psi:M\to M'$ such that this diagram commutes:
  \begin{diagram}[w=3ex,h=3.5ex,tight]
  M    && \rTo^\psi    && M'    \\
  &\luTo    &      & \ruTo  &  \\
  &    & I    & &
  \end{diagram}
  This is the same condition for Saavedra monoids.
\end{blanko}

\section{Gentle multiplicative functors and fair monoidal categories}
%%%%%%%%%%%%%%%%%%%%%%%%%%%%%%%%%%%%%%%%%%%%%%%%%%

\label{Sec:gentle}

In Section~\ref{Sec:monoids} we described a monoid in a monoidal
category as a semi-monoid with a gentle semi-monoid homomorphism from
the unit object.  Clearly the notion of monoidal category should be a
categorification of the notion of monoid, so let us reformulate the
notion of Saavedra unit in terms of {\em gentle multiplicative
functors}.  This viewpoint also leads to the  notion of
fair monoidal category, which is the simplest manifestation of the
higher-categorical concept of fair category, cf.~\cite{Kock:fair}.

\begin{deff}
  We call a functor $\UU\to(\CC,\tensor)$ %, $X \mapsto \OV{X}$
  {\em gentle} if the two induced functors
  $$
    \UU\times\CC  \ \longrightarrow \ \CC\times\CC \ 
    \stackrel{\tensor}{\longrightarrow} \ \CC 
  $$
  $$
    \CC\times\UU \ \longrightarrow \ \CC\times\CC
    \stackrel{\tensor}{\longrightarrow} \ \CC
  $$
  are fully faithful.
%   For the first map this means that the map
%   $$
%   \UU(I,J)\times \CC(X,Y) \to \CC(\stregI,\OV{J}) \times \CC(X,Y) \to 
%   \CC(\stregI X , \OV{J}Y)
%   $$
%   is a bijection for all
%   $I,J$ in $\U$ and all $X,Y$ in $\CC$.
\end{deff}

\begin{lemma}\label{gentle-cancellable}
  Suppose $\UU$ is contractible.  Then $\UU\to\CC$, $U \mapsto \OV{U}$
  is gentle if and only if
  $\OV{U}$ is cancellable for some (and hence for all) $U$ in $\UU$.
\end{lemma}

\begin{dem}
  We consider the left-hand conditions only.
  That $\UU\to\CC$ is gentle means that
  $$
   \UU(I,J)\times \CC(X,Y) \to \CC(\stregI X , \OV{J}Y)
  $$
  is a bijection for all $I,J$ in $\UU$ and all $X,Y$ in $\CC$.  But
  since $\UU$ is contractible the set $\UU(I,J)$ contains precisely
  one element $\psi$, so the bijection is just $\{\psi\}\times\CC(X,Y)
  \to \CC(\stregI X , \OV{J}Y)$, given by tensoring with $\OV{\psi}$
  on the left.  This is to say that $\OV{\psi}$ is tensor cancellable.
  It remains to notice that since $\UU$ is contractible, and since
  isomorphisms between cancellable objects are always tensor
  cancellable, $\OV{\psi}$ is tensor cancellable for all arrows $\psi$
  in $\UU$ if and only if $\stregI$ is cancellable for all objects $I$
  in $\UU$.
%   
%   
%   Suppose $\stregI$ is cancellable in $\CC$.  (Since every object in 
%   $\UU$ is isomorphic to $I$ this implies that all image objects are 
%   cancellable.)  Since all arrows in $\UU$ are invertible, all image 
%   arrows are cancellable in $\CC$; hence we have a bijection $\CC(X,Y)
%   \isopil \CC(\stregI X,\OV{J}Y)$ for all $I,J$ in $\UU$ and all $X,Y$ in 
%   $\CC$.  But since $\UU$ is contractible, there is only one arrow in 
%   the set $\UU(I,J)$, so we also have a bijection $\UU(I,J) \times 
%   \CC(X,Y) \isopil \CC(X,Y)$.  Combining these two bijections we get 
%   the bijection expressing that $\UU\to\CC$ is gentle.
\end{dem}

\begin{cor}\label{*toC}
  A Saavedra unit in $\CC$ is the same thing as a gentle, strong
  multiplicative functor $\eta: * \to \CC$.  
\end{cor}

\begin{dem}
  Call the image object $I$; this object is cancellable if and only if
  $\eta$ is gentle.  Specifying the strong multiplicative
  compatibility on $\eta$ is equivalent to giving the iso\-morphism
  $\alpha: II \to I$.  (Note that since $I$ is cancellable, the
  multiplication map $\alpha$ is automatically associative by
  Lemma~\ref{ass+eq}.)
\end{dem}

Similarly, a morphism of Saavedra units in $\CC$ is the same thing as a
multiplicative natural isomorphism between the corresponding functors
$* \to \CC$.

% \bigskip
% 
% More generally, a natural transformation $u: F\Rightarrow G$ is {\em gentle}
% if for all $I,J\in \UU$ the composite of the big diagram is a 
% bijection of sets:
% \begin{diagram}[w=16ex,h=6ex,tight]
% \UU(I,J)\times\CC(X,Y) & \rTo^{F\times \CC}  & 
% \CC(\stregI,\OV{J})\times\CC(X,Y) &&  \\
% \dTo<{G \times \CC}  &    & \dTo<{(\text{postcomp. } u_J) \times \CC} & \rdTo^\tensor & \\
% \CC(\wtil I,\wtil J)\times\CC(X,Y)  & \rTo_{(\text{precomp. } u_I) 
% \times \CC}  & 
% \CC(\stregI,\wtil J)\times\CC(X, Y) &  & 
% \CC(\stregI X,\OV{J}Y) \\
% & \rdTo_{\tensor} && \rdTo_{\tensor} & \dTo<{\text{postcomp. } u_J Y} \\
% && \CC(\wtil I X,\wtil{J}Y) &\rTo_{\text{precomp. } u_I X}
% & \CC(\stregI X,\wtil J Y)
% \end{diagram}
% It follows that if $u$ is point-wise mono or epi, then it is 
% cancellable.  

\begin{blanko}{Fair monoidal categories.} (Cf.~\cite{Kock:fair}.)
  The idea of fair monoidal category is to take the whole contractible
  category of units, instead of choosing an arbitrary unit in it.  In
  the viewpoint of Corollary~\ref{*toC}, a multiplicative category
  $\CC$ is monoidal when it is equipped with a gentle, strong
  multiplicative functor $* \to \CC$.  By definition, a {\em fair
  monoidal category} is a gentle, strict multiplicative functor $\UU
  \to \CC$ between categories with multiplication $\UU$ and $\CC$,
  with $\UU$ contractible.

  Here $\CC$ is thought of as the `underlying category with
  multiplication', while $\UU$ is thought of as its category of units.
  We want $\UU\to \CC$ to be strict because it should be thought of as
  the forgetful functor from the category of units.
\end{blanko}

\begin{blanko}{Fair monoidal functors.}
  Given fair monoidal categories $(U \to C)$ and $(U' \to C')$, a {\em
  (strong) fair monoidal functor} is a pair of strong multiplicative
  functors $(\phi_U, \phi_C)$ such that the diagram
  \begin{equation}\label{C&U}
  \begin{diagram}[w=6ex,h=4.5ex,tight]
  U & \rTo^{\phi_U}  & U'  \\
  \dTo  &    & \dTo  \\
  C  & \rTo_{\phi_C}  & C'
  \end{diagram}
  \end{equation}
  commutes (as strong multiplicative functors).
  (Note that in \cite{Kock:fair}, only strict fair monoidal functors 
  are considered, i.e.~$\phi_U$ and $\phi_C$ are required 
  to be strict multiplicative functors.)
\end{blanko}

\begin{prop}\label{fair}
  There is an equivalence of categories between the category
  $\kat{MonCat}$ of monoidal categories and strong monoidal functors
  and the category $\kat{FairMonCat}$ of fair monoidal categories and
  (strong) fair monoidal functors.
\end{prop}

\begin{dem}
  We describe first the functor $\kat{MonCat} \to \kat{FairMonCat}$
  which is canonical, where\-as the pseudo-inverse depends on a choice.
  Given a monoidal category $(\CC,\tensor, I)$ then the category
  itself is a strict semi-monoidal category $C$.  The category of units
  $U=U(\CC)$ is another strict semi-monoidal category, which is
  contractible, and there is a strict multiplicative functor $U \to
  C$, which is gentle by \ref{gentle-cancellable}, since each unit is
  cancellable.  Hence $U \to C$ is a fair monoidal category.

  Given a multiplicative functor $(\CC,\tensor,I)\to(\DD,\tensor,J)$,
  $X\mapsto\OV{X}$, the existence of a unit compatibility $\phi_0: J
  \isopil \stregI$ means that $\stregI$ is a unit in $\DD$ (by
  Proposition~\ref{phi-unitmorph}), and by Proposition~\ref{extend}, we can
  lift to a multiplicative functor $U(\CC) \to U(\DD)$, by sending any
  unit $(I,\alpha)$ to $(\stregI,\OV{\alpha})$.  (This is the only
  possible way to lift $\CC\to \DD$ to $U(\CC)\to U(\DD)$.)  By
  Proposition~\ref{extend}, the square~(\ref{C&U}) commutes as strong
  multiplicative functor, so we do have a fair monoidal functor.

  Conversely, suppose we are given a fair monoidal category $(U \to C)$.
  Let $\CC$ denote the category with multiplication $C$.  We must
  provide it with a Saavedra unit.  Just {\em choose} an object $I$ in
  $U$ --- since $U$ is contractible, any choice is as good as any
  other choice.  Now there is induced a canonical unit structure on
  its image $\stregI$ in $C$.  Indeed, the fact that $U$ is
  contractible implies that there is a unique isomorphism $\alpha: II
  \to I$, and the fact that $U \to C$ is gentle means that
  $\stregI$ is cancellable, by \ref{gentle-cancellable}.  Hence
  $\stregI$ is a Saavedra unit in $\CC$.

  Suppose now we are given a diagram (\ref{C&U}) of semi-monoidal
  categories, where $U$ and $U'$ are contractible.  We have already
  chosen $I\in U$ and $I'\in U'$ whose images are units in $C$ and
  $C'$.  The arrows $\alpha: II\to I$ and $\alpha': I'I'\to I'$ are
  given uniquely since $U$ and $U'$ are contractible.  Also since $U'$
  is contractible, there is a unique isomorphism $\phi_0: I' \to \stregI$,
  and it is automatically compatible with $\OV{\alpha}$ and $\alpha'$
  because all compatibilities hold in a contractible category.  The
  image of this compatibility diagram in $\CC'$ is the unit
  compatibility.
  
  It is easy to see that the two constructions are pseudo-inverse to 
  each other, in the sense that going back and forth gives something 
  canonically isomorphic to the starting point, in either direction.
\end{dem}

\begin{blanko}{Monoids in fair monoidal categories.}\label{fair-monoid}
  Monoids in fair monoidal categories have the following appealing
  description.  A fair monoidal category is a gentle multiplicative
  functor $\eta: U \to C$ with $U$ contractible.  A monoid in here is
  a gentle semi-monoid homomorphism $\eta(I) \to M$, where $I$ is a
  semi-monoid in $U$ and $M$ is a semi-monoid in $C$.
\end{blanko}

%%%%%%%%%%%%%%%%%%%%%%%%%%%%%%%%%%%%%%%%%%%%%%%%%%
\section{Lax functoriality}
%%%%%%%%%%%%%%%%%%%%%%%%%%%%%%%%%%%%%%%%%%%%%%%%%%
\label{Sec:lax}

Many of the constructions and arguments of Section~\ref{Sec:functoriality}
admit a lax version, but not everything works as well in the lax 
case: the contractibility of the category of units makes it 
behave better under strong functors than under lax ones.  For example,
when $\phi_0: J \to \stregI$ is not invertible, $\stregI$ is not in 
general a unit itself.

The definition of lax monoidal functor in the \LR{} setting
is well-established: simply 
allow $\phi_2 : \OV{X}\OV{Y} \to \OV{XY}$ and $\phi_0:  J \to \stregI$
to be non-invertible, and otherwise keep 
the conditions of \ref{strong-clas}.  One important motivation for 
considering lax monoidal functors is that monoids are a special case:
a monoid in $\CC$ is essentially the same as a lax monoidal functor
$* \to \CC$ (cf.~B\'enabou~\cite{Benabou:CR1964}).

\begin{blanko}{Lax monoidal functors in the Saavedra-unit setting.}
  A  {\em Saavedra-unit compatibility} for a lax multiplicative functor
  $(\CC,I) \to (\DD,J)$, $X 
  \mapsto \OV{X}$ is a gentle semi-monoid homomorphism
  $\phi_0 : J \to \stregI$, such that
  \begin{equation}\label{epi}
  \stregI \OV{X} \to \OV{IX} \ \text{ and } \ \OV{X}\stregI \to 
  \OV{XI} \ \text{ are epi for all } X .
  \end{equation}
  
  Recall that gentle means that the two composites
  \begin{align}\label{gentle-phi}
  J  \stregI \rTo \stregI  \stregI \rTo \OV{II} \rTo \stregI &
\\
\stregI  J \rTo \stregI  \stregI \rTo \OV{II} \rTo \stregI & \notag
  \end{align}
  are monomorphisms.  (The following proposition implies that in fact those
composites are isomorphisms.)
\end{blanko}

\begin{prop}\label{lax-prop}
  The two notions of lax monoidal functor coincide, under the
  correspondence of Proposition~\ref{bij}.
\end{prop}

The epi condition~(\ref{epi}) may appear a little bit strange, and one might hope
it would be unnecessary (i.e.~automatically satisfied).  Note that the
\LR{} unit compatibility of a lax monoidal functor implies that
$\stregI \OV{X} \to \OV{IX}$ and $\OV{X}\stregI \to \OV{XI}$ are epi
for all $X$.  This follows from the
conditions~(\ref{strong-clas-comp}) since $\lambda_X$ and $\rho_X$ are
isomorphisms.
    
Note that for strong monoidal functors, both condition (\ref{epi})
and (\ref{gentle-phi}) are automatic, so the definition reduces to
\ref{strong-Saa} in this case.

%   Note that $\phi_0\OV{X}$ is split mono for all $X$ in $\CC$,
%   but this has no general implication on $\phi_0$ itself\ldots

\begin{dem*}{Proof of Proposition~\ref{lax-prop}.}
  From \LR{} unit compatibility to Saavedra unit compatibility
  works exactly like in the strong case (Lemma~\ref{funct:clas->Saa}),
  noticing that the composites (\ref{gentle-phi}) are just $\lambda$
  and $\rho$, and hence isomorphisms.  We already observed that 
  (\ref{epi}) is automatic for \LR{} unit compatibilities.
  
  The other direction, starting with a Saavedra-unit compatibility and
  showing that it is also an \LR{} compatibility, is a little bit
  different from the strong case.  First we show that the
  compatibility holds in the special case where the object is $I$.
  This is the content of Lemma~\ref{X=I case}, which does not use the
  epi condition~(\ref{epi}).  Then we use the epi condition to deduce
  the general result from this case (Lemma~\ref{X=X}).
\end{dem*}

% 
% To establish the conjecture in the other direction, the
% first step (independent of the epi condition~(\ref{epi})) 
% is to show that the composite of the last condition is
% precisely $\lambda^J_{\stregIscript}$.

\begin{lemma}\label{X=I case}
  The composite 
  \begin{equation}\label{lambdaJI}
    J  \stregI \rTo \stregI  \stregI \rTo \OV{II} \rTo \stregI
  \end{equation}
  is precisely $\lambda^J_{\stregIscript}$.  And similarly the other
  isomorphism is precisely $\rho^J_{\stregIscript}$.  This does not
  depend on the epi condition~(\ref{epi}).
\end{lemma}

\begin{dem}
  We consider compatibility with $\lambda$; compatibility with $\rho$
  is established analogously.  By construction of $\lambda$, we have
  $J\lambda^J_{\stregIscript} = \beta \stregI$, so to establish the
  assertion of the lemma we tensor the composite with $J$ on the left
  and check that it gives $\beta \stregI$.
  
%   Consider this diagram:
%   \begin{diagram}[w=6ex,h=4.5ex,tight]
%   JJ \stregI && \rTo^{\beta \stregIscript}  && J \stregI \\
%   \dTo<{J\phi_0 \stregIscript}&&&& \dTo^{\phi_0 \stregIscript}  \\
%    J\stregI \stregI& \rTo^{\phi_0 \stregIscript \stregIscript} & 
%   \stregI \stregI \stregI & \rTo^{\OVscript \alpha \stregIscript} &\stregI \stregI  \\
%   \dTo<{J \OVscript{\alpha}}&    & \dTo<{\stregIscript \OVscript{\alpha}} && 
%   \dTo>{\OVscript{\alpha}} \\
% J \stregI & \rTo_{\phi_0 \stregIscript} & \stregI \stregI & \rTo_{\OVscript{\alpha}} 
% & \stregI
%   \end{diagram}
%   The right-hand bottom square is the associativity of the image of a 
%   multiplication map; the left-hand bottom square is just naturality 
%   (elevator business).  The top square is 
%   the compatibility condition relating $\alpha$ and $\beta$.
%   
%   But $\phi_0 \stregI$ followed by $\OV \alpha$ is a monomorphism, so 
%   we can cancel it away and conclude that already the left and the 
%   top maps are equal.  But that is the characterising property of 
%   $\lambda^J_{\stregIscript}$.
  Consider the big diagram~(\ref{bigdiagram}) of the proof in the 
  strong case, but with $X=I$. In this case the 
  right-hand composite and the bottom composite both coincide with 
  (\ref{lambdaJI}), which is a monomorphism by assumption. So we 
  can cancel them away and conclude that the lemma holds in the case 
  $X=I$ as claimed.
\end{dem}

% Now we have shown the proposition in the special case of $X=I$.

\begin{lemma}\label{X=X}
  If $\phi_0 : J \to \stregI$ is a Saavedra unit compatibility for a 
  lax monoidal functor $X\mapsto \OV{X}$, then 
  it is also an \LR{} unit compatibility.
\end{lemma}

\begin{dem}
  By the previous lemma, the \LR{} compatibility holds for the 
  object $I$.  The next step is the case of an object of form $IX$:
  we claim that the composite 
  $$
  \theta: \ 
  J \OV{IX}   \stackrel{\phi_0 \OVscript{IX}}{\rTo}
  \stregI \OV{IX}  \stackrel{\phi_2}{\rTo} 
 \OV{IIX}  \stackrel{\OVscript{\lambda^I_{IX}}}{\rTo}  \OV{IX}
$$
is equal to $\lambda^J_{\OVscript{IX}}$.
  Tensor (\ref{lambdaJI}) with $\OV{X}$ on the right, and apply
  $\phi_2$:
\begin{diagram}[w=7ex,h=5ex,tight]
J \stregI \OV{X} & \rTo^{\phi_0 \stregIscript \OVscript{X}}  & \stregI \stregI \OV{X} 
& \rTo^{\phi_2 \OVscript{X}} & \OV{II} \OV{X} &\rTo^{\OVscript{\alpha} 
\OVscript{X}} & \stregI \OV{X}\\
\dTo<{J \phi_2}&    & \dTo<{\stregIscript \phi_2} && \dTo>{\phi_2} && \dTo>{\phi_2}  \\
J \OV{IX}  & \rTo_{\phi_0 \OVscript{IX}}  & \stregI \OV{IX} & \rTo_{\phi_2} 
& \OV{IIX} & \rTo^{\OVscript{\alpha X}}_{\OVscript{\lambda^I_{IX}}} & \OV{IX}
\end{diagram}
The first square is trivially commutative; the second is 
associativity of $\phi_2$ (cf.~(\ref{phi2-ass})); the third is 
naturality.
The top row is $\lambda^J_{\stregIscript} \OV{X} = 
\lambda^J_{\stregIscript \OVscript{X}}$ (by basic property (2) of $\lambda$).
But we also have the naturality square 
\begin{diagram}[w=6ex,h=4.5ex,tight]
J \stregI \OV{X} & \rTo^{\lambda^J_{\stregIscriptscript \OVscriptscript{X}}}  & \stregI \OV{X}  \\
\dTo<{J \phi_2}  &    & \dTo>{\phi_2}  \\
J \OV{IX}  & \rTo_{\lambda^J_{\OVscriptscript{IX}}}  & \OV{IX}
\end{diagram}
Hence the left-and-bottom ways around in these two diagrams coincide.
Now $\phi_2$ is epi, and tensoring with $J$ on 
the left is an equivalence and hence preserves epimorphisms, so $J\phi_2$
is also epi. So we conclude that $\theta= \lambda^J_{\OVscript{IX}}$.

Now we have proved that the compatibility holds for every object of
form $IX$:
  \begin{diagram}[w=6ex,h=4.5ex,scriptlabels,tight]
  J \OV{IX}       &\rTo^{\lambda^J_{\OVscriptscript{IX}}}    & \OV{IX}    \\
  \dTo<{\phi_0 \OVscript{IX}}    &      & \uTo_{\OVscript{\lambda^I_{IX}}}    \\
  \stregI \OV{IX}&   \rTo_{\phi_2} &  \OV{IIX}
  \end{diagram}

To finish the proof, use the isomorphism $\lambda^I_X : IX \isopil X$
to compare with the diagram we want:
  \begin{diagram}[w=6ex,h=4.5ex,scriptlabels,tight]
  J \OV{X}       &\rTo^{\lambda^J_{\OVscriptscript{X}}}    & \OV{X}    \\
  \dTo<{\phi_0 \OVscript{X}}    &      & \uTo_{\OVscript{\lambda^I_X}}    \\
  \stregI \OV{X}&   \rTo_{\phi_2} &  \OV{IX} .
  \end{diagram}
  
  The compatibility with $\rho$ is established analogously.
\end{dem}
 
\begin{lemma}\label{auto-epi}
  The epi condition~(\ref{epi})
  is automatically satisfied in the following cases:
  \begin{punkt-a}
%     \item $F$ is strongly multiplicative (i.e.~$\phi_2$ are all 
%     isomorphisms).
    
    \item $\phi_0$ is mono and $\phi_{I,X}$ and $\phi_{X,I}$ are mono 
    for all $X$
    
    \item In the domain category, every object is isomorphic to $I$.
%     In particular this is true for $*$, so we get a proof of the 
%     monoid case.
    \end{punkt-a}
\end{lemma}

\begin{dem}
  Re (a):  If both $\phi_0$ and $\phi_{I,X}$ are mono, then the 
  long composite
  $$
  J \OV{X} \to \stregI \OV{X} \to \OV{IX} \to \OV{X}
  $$
  is mono.  But this was all we needed in the proof of 
  Proposition~\ref{compatibilities} (proof of 
  Lemma~\ref{funct:Saa->clas}), so it carries over to the lax 
  case.
  
  Re (b): Lemma~\ref{X=I case} shows that $J\stregI \to 
  \stregI\stregI \to \OV{II} \to \stregI$ is an isomorphism 
  independently of the epi condition~(\ref{epi}).
  Hence $\stregI\stregI \to \OV{II}$ is epi.  But since $X$ is isomorphic
  to $I$, also $\stregI\OV{X} \to \OV{IX}$ is epi.
\end{dem}

\begin{blanko}{Monoids as lax monoidal functors.}
  A monoid in $\CC$ is essentially the same thing as a lax monoidal
  functor $*\to\CC$.  For the \LR{} definition of lax monoidal
  functor this gives the classical notion of monoid, and for
  Saavedra-unit lax monoidal functors this gives the notion of
  Saavedra monoid of \ref{Saavedra-monoid}.  This follows immediately
  since $\OV{*}\OV{*} \to \OV{**}$ is automatically epi, cf.~item (b) 
  of the
  previous lemma.

  Hence as a corollary we get
  Proposition~\ref{monoid-prop}.\label{monoid-proof}
\end{blanko}

\begin{BM}
  For strong monoidal functors we got
  Corollary~\ref{functor-redundance} for free: compatibility with the
  left constraint implies compatibility with the right constraint.
  This result does not follow in the lax case, since we need both
  compatibilities~(\ref{strong-clas-comp}) in order to be able to
  establish the epi condition~(\ref{epi}).  (The conclusion does hold
  in the special situations of \ref{auto-epi}.)  On the other hand, I 
  don't know of a counter example to the statement in the lax case.
\end{BM}

\small

\noindent
\textsc{Address:}  Departament de Matem\`atiques,
Universitat Aut\`onoma de Barcelona, \linebreak 08193 Bellaterra (Barcelona),
Spain.

\noindent
\textsc{Email address:}
  \texttt{kock@mat.uab.cat}

\label{lastpage}

\end{document}